\documentclass[a4paper,12pt]{article}
\usepackage{amsmath}
\usepackage{amssymb}
\usepackage{tabularx}
\usepackage{enumerate}
\usepackage{graphicx}
\usepackage{subfigure}
\usepackage{color}
\usepackage[pagewise]{lineno}
\usepackage{longtable}
\usepackage{epstopdf}
\usepackage{float}
\usepackage{setspace}
\usepackage{mathrsfs}

\usepackage[justification=centering]{caption}
\usepackage{multirow}
\usepackage{cases}
\usepackage[compress]{cite}
\usepackage{lineno}
\usepackage{color}
\usepackage{setspace}
\usepackage{multirow}
 \usepackage{rotating}
\usepackage{geometry}
\newtheorem{thm}{Theorem}[section]
\newtheorem{lem}[thm]{Lemma}
\newtheorem{cor}{Corollary}[section]

\newtheorem{prop}[thm]{Proposition}
\newtheorem{rmk}{Remark}[section]

\newtheorem{defi}{Definition}[section]
\newtheorem{proof}{Proof}

\newcommand{\qed}{\hspace{1em}\mbox{\raisebox{0.65ex}{\fbox{}}}}

\numberwithin{equation}{section}

\newcommand{\be}{\begin{equation}}
\newcommand{\ee}{\end{equation}}
\newcommand\bes{\begin{eqnarray}} \newcommand\ees{\end{eqnarray}}
\newcommand{\bess}{\begin{eqnarray*}}
\newcommand{\eess}{\end{eqnarray*}}

\newcommand{\bpf}{{\bf Proof.\ \ }}
\newcommand{\epf}{\mbox{}\hfill $\Box$}

\begin{document}

\thispagestyle{empty}
\title{On an age-structured model in moving boundaries: The effects of nonlocal diffusion and  harvesting pulse\thanks{ Xu is supported by the Natural Science Foundation of Jiangsu Province, PR China (No. BK20220553);
Santos is supported by the CNPq/Brazil Proc.  $N^{o}$ $311562/2020-5$;
Zhang is supported by the National Natural Science Foundation of China (No. 12101301);
Lin is supported by the National Natural Science Foundation of China (No. 12271470).}}
\date{\empty}
\author{\thanks{Corresponding author. Email: zglin@yzu.edu.cn (Z. Lin).}  \\
{\small $^1$ School of Mathematical Science, Yangzhou University, Yangzhou 225002, China}\\
}

\author{Haiyan Xu$^1$, Carlos Alberto Santos$^2$, Mengyun Zhang$^3$, Zhigui Lin$^1$\thanks{Corresponding author. }  \\
{\small $^1$ School of Mathematical Science, Yangzhou University, Yangzhou 225002, China}\\
{\small $^2$ Department of Mathematics, University of Brasilia, BR-70910900 Brasilia, DF, Brazil}\\
{\small $^3$ School of Applied Mathematics, Nanjing University of Finance Economics, Nanjing 210003, China}
}

\maketitle
\begin{quote}
\noindent
{\bf Abstract} { 
 \small In order to understand how nonlocal diffusion and pulse intervention affect dynamics of species, we focus on an age-structured nonlocal diffusion model in moving and heterogeneous environment, where nonlocal diffusion describes the long range dispersal of species itself and time-periodic harvesting pulse exacting on the adult reflects human intervention. A generalized principal eigenvalue involving harvesting rate used to identify the spreading and vanishing outcomes is firstly defined and the existence of the principal eigenvalue is given under some conditions. Subsequently, properties of the generalized principal eigenvalue and the principal eigenvalue related to harvesting rate and length of habitat interval are analyzed, respectively. The criteria to governing spreading or vanishing of the species are finally investigated, with sufficient conditions for spreading-vanishing established. Our results indicate that complexities can be induced by the internal long rang dispersal and expanding capacities of species, as well as external harvesting intervention of human. Specifically, appropriate harvesting rate and expanding capacities can even change the reciprocal outcomes of species from co-existence to co-extinction.}

\end{quote}

\section{Introduction}
To precisely describe the spatial spreading of species and the front of expanding habitats, several reaction-diffusion mutualistic models of random dispersal with Stefan boundary conditions have been widely researched \cite{LL,ZW}, see also \cite{W,WZ} for competition models and \cite{YMB,WZD} for prey-predator models. The well-posedness of the solution has been established and dynamic behaviors as well as sufficient conditions for spreading-vanishing have been given by a series of well-developed theories since the early work of Du and Lin in \cite{DL}.

Nonlocal diffusion as a long range dispersal can better characterize the movement of species, while local diffusion as adjacent random dispersal processes expressed by reaction-diffusion equations has been widely researched \cite{N,BLS,BC,K}. Some classic models characterized by Lapician local operator have been expanded to nonlocal operators, for instance, a Fisher-KPP equation with nonlocal dispersal in \cite{CD} is a natural extension of the model for a moving boundary in\cite{DL}, which arises significant differences on the spreading-vanishing criteria, and a higher-dimensional Fisher-KPP equation with nonlocal dispersal and a moving boundary is investigated by Du and Ni in \cite{DNN}. See also extensions for a class of two-species mutulistic model \cite{DN,XL}, for models with Lotka-Volterra type competition and prey-predator growth terms \cite{DWZ}, which also give the longtime asymptotic limits of the solution when spreading happens in weak competition and weak predation cases, see also \cite{BL,HN,CLW} and references therein. Besides mathematics, nonlocal diffusion problems have attached much attention in many other fields. A generalized nonlocal Kuramoto-Sivashinsky evolution equation was researched in hydromechanics \cite{DD}. Nonlocal operator method was applied to solve the differential electromagnetic vector wave equations in electric fields \cite{RR}, and nonlocal discrete differential operators as well as a family of weighted p-Laplace operators were introduced to show a simple graph colorization method \cite{LT}. See also \cite{P} for materials science and \cite{ZC} for image recovery.

A nonlocal dispersal operator to replace the local diffusion Laplace term $du_{xx}$ is commonly written by
$$
d(J*u-u)(t,x)=d[\int_{\mathbb{R}}J(x-y)u(t,y)dy-u(t,x)].
$$
Since $u(t,x)$ is the density of species at time $t$ and space $x$, $J(x-y)$ represents the probability distribution function moving from the location $y$ to $x$, then $J(x-y)u(t,y)$ is the rate that is proportional to the probabilities of an individual density $u(t,y)$ jumping from location $y$ to $x$, vice versa. Du et al. in \cite{DLZ,DN,DNNN} established a semi-wave model and discussed different conditions as
$$\int_0^{+\infty} xJ(x)dx<+\infty\,\,\,\mbox{and}\,\,\,\int_0^{+\infty} e^{\lambda x}J(x)dx<+\infty$$
satisfied by a continuous nonnegative even kernel function $J$, that is, spreading speed is finite if and only if $\int_0^{+\infty} xJ(x)dx<+\infty$, otherwise, the accelerate spreading happens and comprehensive results on spreading speed occurs.

Apart from moving boundaries and nonlocal diffusion caused by species itself, the pulses from human intervention also have an effect on the extinction or persistence of species. Mathematical models with periodic pulse interventions such as spraying pesticide and harvesting of crops \cite{TC}, vaccination of susceptible individuals \cite{TM,V} and releasing of natural enemies in the management of invasive species \cite{RN} from continuous to discrete system, will cause considerable technical difficulties for analysis. An impulsive reaction-diffusion model was introduced and its speed and critical threshold were investigated by Lewis and Li in \cite{L}. See also \cite{M} for one-single logistic model with pulses and free boundaries, \cite{XLS} for cooperative reaction-diffusion model with impulsive harvesting in evolving domain, \cite{PL} for a West Nile virus nonlocal model with seasonal succession, and \cite{ZWW  } for competition model with seasonal succession. Taking a comprehensive consideration of these factors, we consider a cooperative nonlocal diffusion model with harvesting pulses imposed on adults in moving and heterogeneous environment, which reads as
\begin{eqnarray}
\left\{
\begin{array}{lll}
u_{1t}=d_1\mathcal{L}_1[u_1]+{b(t)}u_2-[a(t)+m_1(t)]u_1\\[1mm]
\qquad-\alpha_1(t)u_1^2,  & t\in((n\tau)^+,(n+1)\tau],\,x \in (g(t),h(t)),\\[1mm]
u_{2t}=d_2\mathcal{L}_2[u_2]+{a(t)}u_1-m_2(t)u_2-\alpha_2(t)u_2^2, & t\in((n\tau)^+,(n+1)\tau],\,x \in (g(t),h(t)),\\[1mm]
u_1((n\tau)^+,x)=u_1(n\tau,x), & x \in (g(n\tau),h(n\tau)), \\[1mm]
u_2((n\tau)^+,x)=H(u_2(n\tau,x)), & x \in (g(n\tau),h(n\tau)), \\[1mm]
u_1(t,x)=u_2(t,x)=0, & t\in(0,+\infty),\,x \in \{g(t),h(t)\}, \\[1mm]
h'(t)=\sum\limits_{i=1}^{2}\mu_i\int_{g(t)}^{h(t)}\int_{h(t)}^{+\infty}J_i(x-y)u_i(t,x)dydx, &  t\in(n\tau,(n+1)\tau], \\[1mm]
g'(t)=-\sum\limits_{i=1}^{2}\mu_i\int_{g(t)}^{h(t)}\int_{-\infty}^{g(t)}J_i(x-y)u_i(t,x)dydx, &  t\in(n\tau,(n+1)\tau], \\[1mm]
h(0)=-g(0)=h_0, \\[1mm]
u_1(0,x)=u_{1,0}(x), u_2(0,x)=u_{2,0}(x),& x\in[-h_0,h_0],\,n=0,1,2,\dots,\\[1mm]
\end{array} \right.
\label{a01}
\end{eqnarray}
where $u_1(t,x)$ and $u_2(t,x)$ are densities of juveniles and adults, $d_1$ and $d_2$ represent the diffusive rates of juveniles and adults, respectively. $b(t)$ denotes reproduction rate of adults and $a(t)$ is the rate at which juveniles mature into adults. $m_1(t)$ and $m_2(t)$ represent the death rates of juveniles and adults. $\alpha_1(t)$ and $\alpha_2(t)$ denote the competition coefficients of juvenile and adult individuals, respectively. Biologically, $u((n\tau)^+,x)$ can describe the density of survivor mosquito after being killed at every time $t=n\tau$.
Impulsive functions usually take the form $H_1(u)=cu$, $H_2(u)=mu/(a+u)$ and $H_3(u)=ue^{r-bu}$, which are normally called linear function, Beverton-Holt function \cite{B} and Ricker function \cite{R}, respectively. Harvesting function satisfies $H(u)\in C^1([0,+\infty))$ and usually has the properities as following\\

($\mathcal{A}$)\,\,$H(u)>0$ for $u>0$ and $H(0)=0$; $H'(u)>0$ for $u\geq0$. Also, $H(u)/u$ is nonincreasing with respect to $u$ and $0<H(u)/u<1$.

It is clear that different parameters in $H(u)$ will cause different harvesting pulses, as a result, $1-H(u)/u$ is naturally used to characterize the harvesting rate on adults, specially, for $H_1(u)=cu$, $1-c$ is the harvesting rate and $0<1-c<1$. $J_i:\mathbb{R}\rightarrow\mathbb{R}$ with $i=1, 2$ is a continuous nonnegative even kernel function, and has the properties
$$
\mathbf{(J)}: J_i\in C(\mathbb{R})\cap L^{\infty}(\mathbb{R}), J_i(0)>0, \int_{\mathbb{R}}J_i(x)dx=1, \sup_{\mathbb{R}}J_i<\infty, i=1,2,
$$
and nonlocal operator is
\begin{equation}
\mathcal{L}_i[u_i]=\mathcal{L}_i[u_i;g,h](t,x)=\int_{g(t)}^{h(t)}J_i(x-y)u_i(t,y)dy-u_i(t,x), i=1,2.
\label{a001}
\end{equation}
Initial functions $u_{i,0}(x)(i=1,2)$ satisfy
\begin{equation}
u_{i,0}(x)\in C([-h_0,h_0]),\,u_{i,0}(x)>0\,\mbox{in}\,(-h_0,h_0),\,u_{i,0}(-h_0)=u_{i,0}(h_0)=0. \\[1mm]
\label{a02}
\end{equation}

\vspace{3mm}
There have been much research on free boundary problems for mutualistic models. Nevertheless, the introduction of nonlocal dispersal in such an age-structured impulsive system \eqref{a01} in heterogeneous and moving environment is to allow diversity in results and naturally raises analytical difficulties. Whether does the principal eigenvalue that divides dynamic behaviors of species still exist in dual effects of nonlocal dispersal and periodic pulses? What are the new criteria for spreading and vanishing of individuals? Whether or not can external periodic harvesting pulse and internal expanding capacities alter the population state, from persistence to extinction? Moreover, how does the harvesting rate affect or change the dynamical outcomes? We try to address problems raised above. The goal of this paper is to investigate the impact of periodic harvesting pulse and nonlocal diffusion on longtime behavior of species, and finally obtain the spreading-vanishing criteria by overcoming difficulties caused by pulse and nonlocal operator.

The rest of this paper is arranged as follows. In Section 2, the global existence and uniqueness of the solution to problem \eqref{a01} with nonlocal diffusion and harvesting pulse is given. To overcome the difficulty caused by nonlocal operator and harvesting pulse in the periodic eigenvalue problem, we first define a generalized principal eigenvalue involving pulses in Section 3, then sufficient conditions for the existence of the principal eigenvalue are given and the properties of principal eigenvalues are explored. Section 4 mainly focuses on the spreading-vanishing criteria, in which the dynamic behavior of the solution in a fixed domain is firstly investigated, then some sufficient conditions involving pulse on spreading-vanishing are shown. Criteria about expanding capacities governing spreading or vanishing of species are finally obtained.

\section{Global existence and uniqueness}
For convenience, some notations are firstly introduced. For any given $h_0$ and $\tau$, we define

$$\mathbb{H}_{h_0,\tau}:=\{h(t)\in C([0,+\infty))\cap C^1((n\tau,(n+1)\tau]): h(0)=h_0,\,h'(t)>0\},$$
$$\mathbb{G}_{h_0,\tau}:=\{g(t)\in C([0,+\infty))\cap C^1((n\tau,(n+1)\tau]): -g(t)\in \mathbb{H}_{h_0,\tau}\},$$
\begin{equation*}
\begin{array}{ll}
\mathbb{X}_{u_{i,0};\tau}^{g,h}:=&\{(\zeta_1,\zeta_2)(t,x):\zeta_i\in [C((n\tau,(n+1)\tau]\times[g(t),h(t)]), \,\zeta_i\geq 0,\,
\zeta_i(0,x)=u_{i,0}(x)\,\\[2mm]
&\mbox{for}\,x\in[-h_0, h_0]\,\mbox{and}\,\zeta_i(t, g(t))=\zeta_i(t, h(t))=0 \,\mbox{for}\, t\in[0, +\infty),\,i=1,2\}.
\end{array}
\end{equation*}

\begin{thm}\label{th2.1}
Suppose that $\mathbf{(J)}$ holds. Then for any given $h_0$, $u_{i,0}$ satisfying \eqref{a02} and any given $\tau>0$, problem \eqref{a01} admits a unique solution $(u_1,u_2;(g,h))$ defined for all $t\in(0,+\infty)$. Moreover,
$$
(u_1(t,x),u_2(t,x))\in [\mathbb{X}_{u_{i,0};\tau}^{g,h}]^2\,\,\mbox{and}\,\,(g(t),h(t))\in \mathbb{G}_{h_0,\tau} \times \mathbb{H}_{h_0,\tau}.
$$
\end{thm}
\bpf
(i) For the interval $t\in(0^+,\tau]$, we regard $(u_1(0^+,x),u_2(0^+,x))$ as the initial value of solution $(u_1(t,x),u_2(t,x))$ to problem \eqref{a01}. Since $u_{i,0}(x)\in C([-h_0,h_0])$ and $H\in C^1([0,+\infty))$, we derive that new initial value satisfies that $u_1(0^+,x)=u_1(0,x)\in C([-h_0,h_0])$ and $u_2(0^+,x)=H(u_{2,0}(x))\in C([-h_0,h_0])$. Recalling the existence and uniqueness of the solution without pulse in \cite{CD} and Lemma \ref{m2} later, we declare that the solution $(u_1,u_2;(g,h))$ to problem \eqref{a01} exists and is unique in $t\in(0^+,\tau]$.  Furthermore, $u_i\in C((0,\tau]\times[g(t),h(t)])$ for $i=1,2$ and $(g(t),h(t))\in [C([0,\tau])\cap C^1((0,\tau])]^2$.

(ii) For the interval $t\in(\tau^+,2\tau]$, by the same procedure as (i), the new initial value satisfies $u_1(\tau^+,x)=u_1(\tau,x)\in C([g(\tau),h(\tau)])$ and $u_2(\tau^+,x)=H(u_2(\tau,x))\in C([g(\tau),h(\tau)])$. Owing to Lemma \ref{m2}, one easily checks that the solution $(u_1,u_2;(g,h))$ to problem \eqref{a01} exists and is unique in $t\in(\tau^+,2\tau]$. Also, $(u_1(t,x),u_2(t,x))\in [C((\tau,2\tau]\times[g(t),h(t)])]^2$ and $(g(t),h(t))\in [C([\tau,2\tau])\cap C^1((\tau,2\tau])]^2$.

(iii) The local existence and uniqueness of the solution can be derived by the same process in interval $t\in(2\tau^+,3\tau]$, $t\in(3\tau^+,4\tau]$,
\dots, and step by step, we then find a maximal time interval $[0,T_{\max})$ with $T_{\max}:=n_0\tau+\tau_0$, $0\leq\tau_0<\tau$ and positive integer $n_0$ by Zorn's lemma such that problem \eqref{a01} admits a unique solution in $[0,T_{\max})$.

(iv) We now claim that $T_{\max}=+\infty$. The following estimates of $(u_1,u_2;(g,h))$ in Lemma \ref{m2}, which together with the standard continuous extension method, yield the global existence and uniqueness of the solution to problem \eqref{a01}.

\epf

\begin{lem}
Assume that $\mathbf{(J)}$ holds. For any given nonnegative integer $n_0$ and $0<\tau_1\leq \tau$, if $(u_1,u_2;(g,h))$ is a solution to problem \eqref{a01} defined for $t\in(0,T]$ with $T:=n_0\tau+\tau_1$. Then we acquire
$$
0<u_1(t,x),u_2(t,x)\leq\max\{\frac{b^M}{\alpha^m _1},\frac{a^M}{\alpha^m_2},\parallel u_{1,0}\parallel_\infty,\parallel u_{2,0}\parallel_\infty\}:=A
$$
for $t\in(0,T]$ and $x\in(g(t),h(t))$.
\label{m2}
\end{lem}
\bpf
Denote $(\bar u_1,\bar u_2)=(A,A)$, since $\mathbf{(J)}$ holds, careful calculations yield
$$\begin{array}{llllll}
&&\bar u_{1t}-d_1\mathcal{L}_1[\bar u_1]-{b(t)}\bar u_2+(a(t)+m_1(t))\bar u_1-\alpha_1(t)\bar u_1^2\\[1mm]
&\geq&-d_1A(\int_{g(t)}^{h(t)}J_1(x-y)dy-1)-{b(t)}A-\alpha_1(t)\bar A^2\\[1mm]
&\geq&A(-b^M+\alpha_1^mA)\geq0,\\[1mm]
\end{array}$$
similarly, $\bar u_{2t}-d_2\mathcal{L}_2[\bar u_2]-{a(t)}\bar u_1+m_2(t)\bar u_2-\alpha_2(t)\bar u_2^2>0$ for $t\in(0,T]$ and $x\in(g(t),h(t))$.

If $n_0=0$, we get $t\in (0,T]\subseteq (0,\tau]$.  Since $(\bar u_1(0,x),\bar u_2(0,x))\geq (u_{1,0},u_{2,0})$, a comparison principle to conclude $(0,0)<(u_1(t,x),u_2(t,x))\leq (A,A)$ for $t\in (0,\tau_1]$ and $x\in(g(t),h(t))$. If $n_0=1$, we obtain $t\in(0,\tau+\tau_1]$. Since the interval $t\in (0,\tau_1]$ is discussed above, we here fix $t\in ((\tau_1)^+,\tau+\tau_1]$. It then follows from ($\mathcal{A}$) that
$$
\bar u_1((\tau_1)^+,x)=\bar u_1(\tau_1,x),\,\, \bar u_2((\tau_1)^+,x)=A>H(A)=H(\bar u_2(\tau_1,x)),x\in(g(n\tau_1),h(n\tau_1)).
$$
Step by step, $(0,0)<(u_1(t,x),u_2(t,x))\leq (A,A)$ for $t\in (0,T]$ and $x\in(g(t),h(t))$ can be obtained.
\epf
\section{The generalized eigenvalue problem}
To understand the long-time behavior of the solution to problem \eqref{a01}, we first consider the corresponding problem of \eqref{a01} in a fixed interval $[L_1,L_2]$, and the following nonlocal time-periodic eigenvalue problem in a fixed boundary
\begin{eqnarray}
\left\{
\begin{array}{lll}
\phi_t-d_1\mathcal{\tilde{L}}_1[\phi]=b(t)\psi-[a(t)+m_1(t)]\phi+\lambda \phi,& t\in(0^+,\tau],\,\,x \in (L_1,L_2),  \\[1mm]
\psi_t-d_2\mathcal{\tilde{L}}_2[\psi]=a(t)\phi-m_2(t)\psi+\lambda \psi, & t\in(0^+,\tau],\,\,x \in (L_1,L_2),  \\[1mm]
\phi(0^+,x)=\phi(0,x),& x\in (L_1,L_2),\\[1mm]
\psi(0^+,x)=H'(0)\psi(0,x),& x\in (L_1,L_2),\\[1mm]
\phi(0,x)=\phi(\tau,x),\psi(0,x)=\psi(\tau,x), & x\in [L_1,L_2]
\end{array} \right.
\label{c02}
\end{eqnarray}
is firstly introduced, where $\mathcal{\tilde{L}}_i[u]$ is defined in \eqref{a001} with $[g(t),h(t)]$ replaced by $[L_1,L_2]$, satisfying

$$\mathcal{\tilde{L}}_i[u]:=\mathcal{\tilde{L}}_i[u;L_1,L_2](t,x)=\int_{L_1}^{L_2}J_i(x-y)u(t,y)dy-u(t,x), i=1,2.$$

\vspace{2mm}

For the coupled problem \eqref{c02} with harvesting pulse, the generalized principal eigenvalues can be defined by $\overline{\lambda}$ and $\underline{\lambda}$:
$$\begin{array}{llllll}
&&\overline{\lambda}((L_1,L_2),H'(0)):=\inf\{\lambda\in \mathbb{R}\,|\, \exists \, \phi, \psi\in C^{1, 0}((0,\tau]\times [L_1, L_2]), \phi, \psi>0\,\mbox{and}\, \phi, \psi\, \mbox{are}\\[1mm]
&&\, \tau-\mbox{periodic}\,\mbox{so as \eqref{c01} hold}\},
\end{array}$$
where \begin{eqnarray}
\left\{
\begin{array}{lll}
\phi_t-d_1\mathcal{\tilde{L}}_1[\phi]\leq b(t)\psi-[a(t)+m_1(t)]\phi+\lambda \phi,& t\in(0^+,\tau],\,\,x \in (L_1,L_2),  \\[1mm]
\psi_t-d_2\mathcal{\tilde{L}}_2[\psi]\leq a(t)\phi-m_2(t)\psi+\lambda \psi, & t\in(0^+,\tau],\,\,x \in (L_1,L_2),  \\[1mm]
\phi(0^+,x)=\phi(0,x),& x\in (L_1,L_2),\\[1mm]
\psi(0^+,x)=H'(0)\psi(0,x),& x\in (L_1,L_2),\\[1mm]
\phi(0,x)=\phi(\tau,x),\psi(0,x)=\psi(\tau,x), & x\in [L_1,L_2],
\label{c01}
\end{array} \right.
\end{eqnarray}
and
$$\begin{array}{llllll}
&&\underline{\lambda}((L_1,L_2),H'(0)):=\sup\{\lambda\in \mathbb{R}\,|\, \exists \, \phi, \psi\in C^{1, 0}((0,\tau]\times [L_1, L_2]), \phi, \psi>0\, \mbox{and}\, \phi, \psi\, \mbox{are}\,\\[1mm]
&&\tau-\mbox{periodic if inequalities in \eqref{c01} are all reversed}.\}\\[1mm]
\end{array}$$
If we can find $(\lambda^*,\phi,\psi)$ with positive function pair $(\phi,\psi)$ such that the equalities in \eqref{c01} hold, then $\lambda^*$ is called a principal eigenvalue of problem \eqref{c02}.

\vspace{3mm}

Specially, suppose $J_1(x)=J_2(x)(:=J(x))$ and all coefficients are constant, then periodic eigenvalue problem \eqref{c02} can be transformed into
\begin{eqnarray}
\left\{
\begin{array}{lll}
\phi_t-d_1\mathcal{\tilde{L}}[\phi]=b\psi-(a+m_1)\phi+\lambda \phi,& t\in(0^+,\tau],\,\,x \in (L_1,L_2),  \\[1mm]
\psi_t-d_2\mathcal{\tilde{L}}[\psi]=a\phi-m_2\psi+\lambda \psi, & t\in(0^+,\tau],\,\,x \in (L_1,L_2),  \\[1mm]
\phi(0^+,x)=\phi(0,x),& x\in (L_1,L_2),\\[1mm]
\psi(0^+,x)=H'(0)\psi(0,x),& x\in (L_1,L_2),\\[1mm]
\phi(0,x)=\phi(\tau,x),\psi(0,x)=\psi(\tau,x), & x\in [L_1,L_2].
\end{array} \right.
\label{c07}
\end{eqnarray}
Moreover, we declare that problem \eqref{c07} can be transformed into two eigenvalue problems since $J_1=J_2$, one of which is time-independent, and the other is space-independent. In fact, we first consider a time-independent eigenvalue problem
$$\mathcal{\tilde{L}}[u]:=\mathcal{\tilde{L}}[u(x);L_1,L_2]=\int_{L_1}^{L_2}J(x-y)u(y)dy-u(x)=\lambda_0u(x),\,L_1\leq x\leq L_2.$$
It follows from Cao et al. in \cite[Proposition 3.4 with $a_0=0$]{CD} that the principal eigen-pair $(\lambda_0,\Psi(x))$ of operator $\mathcal{\tilde{L}}$ exists with positive eigenvalue function $\Psi(x)\in C([L_1,L_2])$. By the method of separation variables with $\phi(t,x)=\alpha(t)\Psi(x)$ and $\psi(t,x)=\beta(t)\Psi(x)$, the eigenvalue problem \eqref{c02}(problem \eqref{c07}) can be written by a spatially independent eigenvalue problem
\begin{eqnarray}
\left\{
\begin{array}{lll}
\alpha'(t)=b\beta(t)-(a+m_1-d_1\lambda_0)\alpha(t)+\lambda\alpha(t),& t\in(0^+,\tau],\\[1mm]
\beta'(t)=a\alpha(t)-(m_2-d_2\lambda_0)\beta(t)+\lambda\beta(t),& t\in(0^+,\tau],\\[1mm]
\alpha(0)=\alpha(\tau),\beta(0)=\beta(\tau),\\[1mm]
\alpha(0^+)=\alpha(0),\beta(0^+)=H'(0)\beta(0).\\[1mm]
\end{array} \right.
\label{c03}
\end{eqnarray}
We derive from the first two equations of problem \eqref{c03} that
\begin{equation*}
\left(
\begin{array}{c}
\alpha'(t)\\
\beta'(t)\\
\end{array}
\right)
=
\left(
\begin{array}{cc}
-a-m_1+d_1\lambda_0+\lambda & b\\
a & -m_2+d_2\lambda_0+\lambda\\
\end{array}
\right)
\left(
\begin{array}{c}
\alpha(t)\\
\beta(t)\\
\end{array}
\right)
\triangleq M
\left(
\begin{array}{c}
\alpha(t)\\
\beta(t)\\
\end{array}
\right).
\end{equation*}
The corresponding characteristic equation is $|M-\mu E|=0$, and a direct calculation yields
\begin{equation}
\mu_{1,2}=\lambda+\frac{-a-m_1-m_2+(d_1+d_2)\lambda_0\underline{+}\sqrt{(a+m_1-m_2-d_1\lambda_0+d_2\lambda_0)^2+4ab}}{2}\triangleq \lambda+c_{1,2}.
\label{bbb}
\end{equation}
Without loss of generality, we assume that $c_1>c_2$. Careful calculations imply that
$$a+m_1-d_1\lambda_0+c_1=-(m_2-d_2\lambda_0+c_2)>0.$$

Let $(k_{11},k_{12})$ and $(k_{21},k_{22})$ be the linearly independent of eigenvectors related to eigenvalues $\mu_1$ and $\mu_2$, which yield
\begin{equation*}
\left(
\begin{array}{cc}
k_{i1} & k_{i2}\\
\end{array}
\right)
\left(
\begin{array}{cc}
-a-m_1+d_1\lambda_0+\lambda-\mu_i & b\\
a & -m_2+d_2\lambda_0+\lambda-\mu_i\\
\end{array}
\right)
=
\left(
\begin{array}{cc}
0 & 0\\
\end{array}
\right)
\end{equation*}
for $i=1,2$. Therefore
$$
(k_{11},k_{12})=(a,a+m_1-d_1\lambda_0-\lambda+\mu_1)=(a,a+m_1-d_1\lambda_0+c_1)
$$
and
$$
(k_{21},k_{22})=(m_2-d_2\lambda_0-\lambda+\mu_2,b)=(m_2-d_2\lambda_0+c_2,b).
$$

In the following, we consider the corresponding algebraic equations
\begin{equation*}
\left(
\begin{array}{cc}
a & a+m_1-d_1\lambda_0+c_1\\
m_2-d_2\lambda_0+c_2 & b\\
\end{array}
\right)
\left(
\begin{array}{c}
\alpha(t)\\
\beta(t)\\
\end{array}
\right)
=
\left(
\begin{array}{c}
e^{\mu_1t}\\
me^{\mu_2t}\\
\end{array}
\right),
\end{equation*}
and through careful calculations we can obtain
$$
(\alpha(t),\beta(t))=
(\frac{be^{\mu_1t}-(a+m_1-d_1\lambda_0+c_1)me^{\mu_2t}}{C},\frac{-(m_2-d_2\lambda_0+c_2)e^{\mu_1t}+ame^{\mu_2t}}{C}),
$$
where
$$C=ab-(a+m_1-d_1\lambda_0+c_1)(m_2-d_2\lambda_0+c_2)>0.$$
Direct calculations yield $\alpha(0)+k_1\beta(0)=k_2$ with $k_1=(a+m_1-d_1\lambda_0+c_1)/a$ and $k_2=b/C-(a+m_1-d_1\lambda_0+c_1)(m_2-d_2\lambda_0+c_2)/(aC)$, which are independent of $m$.

In the following we denote $\Lambda=e^{\mu_1\tau}$ for simplicity. Since $\alpha(0)=\alpha(\tau)$ and $\beta(0^+)=H'(0)\beta(\tau)$, we get
\begin{eqnarray}
\left\{
\begin{array}{lll}
b-(a+m_1-d_1\lambda_0+c_1)m=b\Lambda-(a+m_1-d_1\lambda_0+c_1)e^{(c_2-c_1)\tau}m\Lambda, \\[1mm]
am+(a+m_1-d_1\lambda_0+c_1)=H'(0)ae^{(c_2-c_1)\tau}m\Lambda+H'(0)(a+m_1-d_1\lambda_0+c_1)\Lambda.\\[1mm]
\end{array} \right.
\label{c04}
\end{eqnarray}
For abbreviation, we further denote
$$A_{11}=b, \,A_{12}=a+m_1-d_1\lambda_0+c_1, \,A_{13}=(a+m_1-d_1\lambda_0+c_1)e^{(c_2-c_1)\tau},$$
$$A_{21}=a,\, A_{22}=H'(0)ae^{(c_2-c_1)\tau},\,A_{23}=H'(0)(a+m_1-d_1\lambda_0+c_1),$$
then \eqref{c04} can be written by
\begin{eqnarray}
\left\{
\begin{array}{lll}
A_{11}-A_{12}m=A_{11}\Lambda-A_{13}m\Lambda, \\[1mm]
A_{21}m+A_{12}=A_{22}m\Lambda+A_{23}\Lambda.\\[1mm]
\end{array} \right.
\label{c05}
\end{eqnarray}
The explicit solution $(m,\Lambda)$ to problem \eqref{c05} is not easy to given, so in the following we will develop the existence of solution by image method. It is clear that \eqref{c05} can be transformed into
\begin{eqnarray}
\left\{
\begin{array}{lll}
\Lambda=(A_{11}-A_{12}m)/(A_{11}-A_{13}m), \\[1mm]
\Lambda=(A_{21}m+A_{12})/(A_{22}m+A_{23}).\\[1mm]
\end{array} \right.
\label{c06}
\end{eqnarray}

\begin{figure}
\centering
 {
\includegraphics[width=0.45\textwidth]{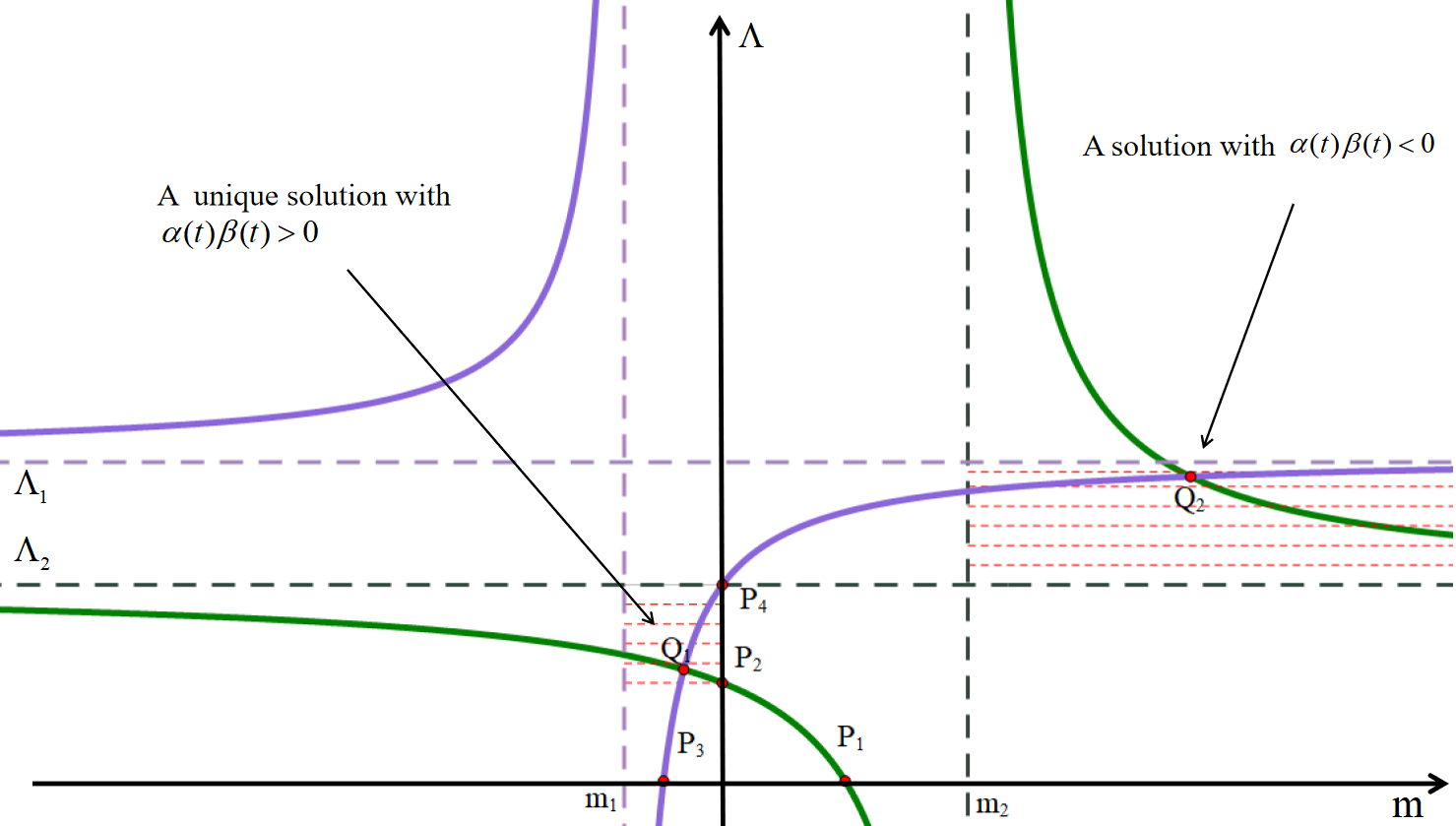}
}
\caption{\scriptsize The solution curve to problem $\eqref{c06}$. The green line represents the solution curve of the first equation in $\eqref{c06}$, which go through fixed points $P_1(\frac{b}{a+m_1-d_1\lambda_0+c_1},0)$ and $P_2(0,1)$. $m_1=-(a+m_1-d_1\lambda_0+c_1)/(ae^{(c_2-c_1)\tau})$ and $\Lambda_1=1/(H'(0)e^{(c_2-c_1)\tau})$ are two asymptotes. Also, $\Lambda$ is strictly decreasing with respect to $m$ in interval $(-\infty,\frac{b}{(a+m_1-d_1\lambda_0+c_1)e^{(c_2-c_1)\tau}})\bigcup((a+m_1-d_1\lambda_0+c_1)e^{(c_2-c_1)\tau},\infty)$ and  $\lim\limits_{m\to\infty}\Lambda=1/e^{(c_2-c_1)\tau}(>1)$. Similarly, mathematical analysis of the second equation in $\eqref{c06}$ shows that points $P_3(-\frac{a+m_1-d_1\lambda_0+c_1}{a},0)$, $P_4(0,\frac{1}{e^{(c_2-c_1)\tau}})$ are in the blue line and  $m_2=b/[(a+m_1-d_1\lambda_0+c_1)e^{(c_2-c_1)\tau}]$, $\Lambda_2=1/(e^{(c_2-c_1)\tau})$ are two asymptotes. $\Lambda$ is strictly increasing with respect to $m$ in interval $(-\infty,\frac{-(a+m_1-d_1\lambda_0+c_1)}{ae^{(c_2-c_1)\tau}})\bigcup(\frac{-(a+m_1-d_1\lambda_0+c_1)}{ae^{(c_2-c_1)\tau}},\infty)$ and  $\lim\limits_{m\to\infty}\Lambda=1/H'(0)e^{(c_2-c_1)\tau}(>1)$. Clearly, $Q_1$ and $Q_2$ are two intersection points, which represents two solution to problem \eqref{c06}.
}
\end{figure}
It can be seen from Fig. 1 that points $Q_1(m_1,\Lambda_1)\in D_1:=(\frac{-(a+m_1-d_1\lambda_0+c_1)}{ae^{(c_2-c_1)\tau}},0)\times(1,\frac{1}{e^{(c_2-c_1)\tau}})$ and $Q_2(m_2,\Lambda_2)\in D_2:=(\frac{b}{(a+m_1-d_1\lambda_0+c_1)e^{(c_2-c_1)\tau}},\infty)\times(\frac{1}{e^{(c_2-c_1)\tau}},\frac{1}{H'(0)e^{(c_2-c_1)\tau}})$. Since $\alpha(t)\beta(t)<0$ in $D_2$, problem \eqref{c04} admits a unique solution $Q_1(m_1,\Lambda_1)\in D_1$ with $\alpha(t)\beta(t)>0$, then eigenvalue problem \eqref{c03} has a unique solution with $\alpha(0)+k_1\beta(0)=k_2$ and $(\alpha(t),\beta(t))>(0,0)$ for $t\in(0^+,\tau]$, which indicates problem \eqref{c07} admits a unique solution $(\lambda^*,\phi,\psi)$ with $\phi,\psi>0$.

\vspace{3mm}

\begin {rmk}\label{rm3.1}
Suppose that $J_1(x)=J_2(x)(:=J(x))$.  If all coefficients are constant, it follows from the analysis in Fig. 1 that the principal eigen-pair $(\lambda^*((L_1,L_2),H'(0)),\phi(t,x),\psi(t,x))$ of periodic problem \eqref{c02} (i.e. problem \eqref{c07}) exists with $\phi(t,x),\psi(t,x)>0$ in $t\in(0^+,\tau]$ and $x\in[L_1,L_2]$. If coefficients are time-dependent, it follows from \cite{BZ} that the principal eigenvalue $\lambda^*((L_1,L_2),H'(0))$ in \eqref{c03} exists with positive eigenfunctions. Otherwise, if $J_1(x)\neq J_2(x)$, then eigenvalue problem \eqref{c02} only admits generalized eigenvalues and we denote them by $\overline{\lambda}((L_1,L_2),H'(0))$ and $\underline{\lambda}((L_1,L_2),H'(0))$, respectively.
\end {rmk}

Next, we give a comparison principle involving the generalized principal eigenvalue, which will be used in the sequel.
\begin{lem} \label{le3.1}Assume that $\mathbf{(J)}$ holds.

$(i)$ If $J_1\neq J_2$ and all coefficients in \eqref{c02} are constant. Then for any given $\overline \lambda_1((L_1,L_2),H'(0))$ and functions $\phi(t,x),\psi(t,x)\in C^{1, 0}((0,\tau]\times [L_1, L_2])$ with $\phi,\psi\geq,\not\equiv0$ satisfying
\begin{eqnarray}
\left\{
\begin{array}{lll}
\phi_t-d_1\mathcal{\tilde{L}}_1[\phi]\leq b\psi-(a+m_1)\phi+\overline{\lambda}_1 \phi,& t\in(0^+,\tau],\,\,x \in (L_1,L_2),  \\[1mm]
\psi_t-d_2\mathcal{\tilde{L}}_2[\psi]\leq a\phi-m_2\psi+\overline{\lambda}_1 \psi, & t\in(0^+,\tau],\,\,x \in (L_1,L_2),  \\[1mm]
\phi(0^+,x)=\phi(0,x),& x\in (L_1,L_2),\\[1mm]
\psi(0^+,x)=H'(0)\psi(0,x),& x\in (L_1,L_2),\\[1mm]
\phi(0,x)=\phi(\tau,x),\psi(0,x)=\psi(\tau,x), & x\in [L_1,L_2],
\end{array} \right.
\label{c11}
\end{eqnarray}
and $(\underline{\lambda}_1,\phi_1,\psi_1)$ satisfying
\begin{eqnarray}
\left\{
\begin{array}{lll}
\phi_{1t}-d_1\mathcal{\tilde{L}}_1[\phi_1]\geq b\psi_1-(a+m_1)\phi_1+\underline{\lambda}_1\phi_1,& t\in(0^+,\tau],\,\,x \in (L_1,L_2),  \\[1mm]
\psi_{1t}-d_2\mathcal{\tilde{L}}_2[\psi_1]\geq a\phi_1-m_2\psi_1+\underline{\lambda}_1\psi_1, & t\in(0^+,\tau],\,\,x \in (L_1,L_2),  \\[1mm]
\phi_1(0^+,x)=\phi_1(0,x),& x\in (L_1,L_2),\\[1mm]
\psi_1(0^+,x)=H'(0)\psi_1(0,x),& x\in (L_1,L_2),\\[1mm]
\phi_1(0,x)=\phi_1(\tau,x),\psi_1(0,x)=\psi_1(\tau,x), & x\in [L_1,L_2],
\end{array} \right.
\label{c12}
\end{eqnarray}
we obtain $\underline{\lambda}_1\leq\overline{\lambda}_1$.

$(ii)$ If $J_1=J_2$. Then for any given $(\overline{\Lambda},\Phi,\Psi)$ with $(\Phi(t,x),\Psi(t,x))\in ([0,\tau]\times[L_1,L_2])^2$ and $(\Phi,\Psi)\geq,\not\equiv(0,0)$, satisfying
\begin{eqnarray}
\left\{
\begin{array}{lll}
\Phi_t-d_1\mathcal{\tilde{L}}[\Phi]\leq b(t)\Psi-[a(t)+m_1(t)]\Phi+\overline{\Lambda}\Phi,& t\in(0^+,\tau],\,\,x \in (L_1,L_2),  \\[1mm]
\Psi_t-d_2\mathcal{\tilde{L}}[\Psi]\leq a(t)\Phi-m_2(t)\Psi+\overline{\Lambda}\Psi, & t\in(0^+,\tau],\,\,x \in (L_1,L_2),  \\[1mm]
\Phi(0^+,x)=\Phi(0,x),& x\in (L_1,L_2),\\[1mm]
\Psi(0^+,x)=H'(0)\Psi(0,x),& x\in (L_1,L_2),\\[1mm]
\Phi(0,x)=\Phi(\tau,x),\Psi(0,x)=\Psi(\tau,x), & x\in [L_1,L_2],
\end{array} \right.
\label{c13}
\end{eqnarray}
and $(\underline{\Lambda},\Phi_1,\Psi_1)$ satisfying
\begin{eqnarray}
\left\{
\begin{array}{lll}
\Phi_{1,t}-d_1\mathcal{\tilde{L}}[\Phi_1]\geq b(t)\Psi_1-[a(t)+m_1(t)]\Phi_1+\underline{\Lambda}\Phi_1,& t\in(0^+,\tau],\,\,x \in (L_1,L_2),  \\[1mm]
\Psi_{1,t}-d_2\mathcal{\tilde{L}}[\Psi_1]\geq a(t)\Phi_1-m_2(t)\Psi_1+\underline{\Lambda}\Psi_1, & t\in(0^+,\tau],\,\,x \in (L_1,L_2),  \\[1mm]
\Phi_1(0^+,x)=\Phi_1(0,x),& x\in (L_1,L_2),\\[1mm]
\Psi_1(0^+,x)=H'(0)\Psi_1(0,x),& x\in (L_1,L_2),\\[1mm]
\Phi_1(0,x)=\Phi_1(\tau,x),\Psi_1(0,x)=\Psi_1(\tau,x), & x\in [L_1,L_2],
\end{array} \right.
\label{c133}
\end{eqnarray}
we have $\underline{\Lambda}\leq \overline{\Lambda}$, and the equality holds provided that $(\overline{\Lambda},\Phi,\Psi)$ and $(\underline{\Lambda},\Phi_1,\Psi_1)$ are the principal eigen-pair of eigenvalue problem \eqref{c02}.
\end{lem}
\bpf
Before proving, we first define $$\langle\mu,\nu\rangle:=\int_{L_1}^{L_2}\mu(t,x)\nu(t,x)dx.$$

(i) If $J_1\neq J_2$ and all coefficients in \eqref{c01} are constant. Let $W(t,x)=\phi(\tau-t,x)$ and $Z(t,x)=\psi(\tau-t,x)$ in \eqref{c11}, which satisfy
\begin{eqnarray}
\left\{
\begin{array}{lll}
-W_t-d_1\mathcal{\tilde{L}}_1[W]\leq bZ-(a+m_1)W+\overline{\lambda}_1 W,& t\in[0,\tau^-),\,\,x \in (L_1,L_2),  \\[1mm]
-Z_t-d_2\mathcal{\tilde{L}}_2[Z]\leq aW-m_2Z+\overline{\lambda}_1 Z, & t\in[0,\tau^-),\,\,x \in (L_1,L_2),  \\[1mm]
W(\tau^-,x)=W(\tau,x),& x\in (L_1,L_2),\\[1mm]
Z(\tau^-,x)=H'(0)Z(\tau,x),& x\in (L_1,L_2),\\[1mm]
W(0,x)=W(\tau,x),Z(0,x)=Z(\tau,x), & x\in [L_1,L_2],
\end{array} \right.
\label{c166}
\end{eqnarray}
where $W(\tau^-,x)=\lim\limits_{t\to\tau-0}{W(t,x)}$ and $Z(\tau^-,x)=\lim\limits_{t\to\tau-0}{Z(t,x)}$ denote the left limits of $W(t,x)$ and $Z(t,x)$ at $t=\tau$, respectively.
Firstly, multiplying the first equation in \eqref{c12} by $W$ and the first equation in \eqref{c166} by $\phi_1$, then integrating both sides of these two equations in $x$ over $[L_1,L_2]$ give that
$$
\langle\phi_{1t},W\rangle-d_1\langle\mathcal{\tilde{L}}_1[\phi_1],W\rangle\geq b\langle \psi_1,W\rangle-(a+m_1)\langle\phi_1,W\rangle+\underline{\lambda}_1 \langle\phi_1,W\rangle,\,t\in(0,\tau),
$$
$$
\langle-W_t,\phi_1\rangle-d_1\langle\mathcal{\tilde{L}}_1[W],\phi_1\rangle\leq b\langle Z,\phi_1\rangle-(a+m_1)\langle W,\phi_1\rangle+\overline{\lambda}_1\langle W,\phi_1\rangle,\,t\in(0,\tau).
$$

In the following, since $\langle\mathcal{\tilde{L}}_1[\phi_1],W\rangle=\langle\mathcal{\tilde{L}}_1[W],\phi_1\rangle$, we obtain by abstracting these two equations and integrating in $t$ over $(0^+,\tau^-)$ that
$$
\int_{o^+}^{\tau^-}\langle\phi_{1t},W\rangle+\langle W_t,\phi_1\rangle dt\geq
\int_{0^+}^{\tau^-}b(\langle \psi_1,W\rangle-\langle Z,\phi_1\rangle)dt+(\underline{\lambda}_1-\overline{\lambda}_1)\int_{0^+}^{\tau^-}\langle\phi_1,W\rangle dt.
$$
Since
$$\begin{array}{llllll}
&&\int_{0^+}^{\tau^-}\langle\phi_{1t},W\rangle+\langle W_t,\phi_1\rangle dt\\[1mm]
&=&\int_{L_1}^{L_2}(\phi_1(\tau^-,x)W(\tau^-,x)-\phi_1(0^+,x)W(0^+,x))dt&\\[1mm]
&=&\int_{L_1}^{L_2}(\phi_1(0^+,x)W(\tau^-,x)-\phi_1(0^+,x)W(\tau^-,x))&\\[1mm]
&=&0,
\end{array}$$
thus
\begin{equation}
0\geq b\int_{0^+}^{\tau^-}(\langle \psi_1,W\rangle-\langle Z,\phi_1\rangle)dt+(\underline{\lambda}_1-\overline{\lambda}_1)\int_{0^+}^{\tau^-}\langle\phi_1,W\rangle dt.
\label{c17}
\end{equation}
Similarly, repeating the above process for the second equations in \eqref{c12} and \eqref{c166}, it then follows from $\langle\mathcal{\tilde{L}}_2[\psi_1],Z\rangle=\langle\mathcal{\tilde{L}}_2[Z],\psi_1\rangle$ and
$$\begin{array}{llllll}
&&\int_{0^+}^{\tau^-}\langle\psi_{1t},Z\rangle+\langle Z_t,\psi_1\rangle dt\\[1mm]
&=&\int_{L_1}^{L_2}(\psi_1(\tau^-,x)Z(\tau^-,x)-\psi_1(0^+,x)Z(0^+,x))dt&\\[1mm]
&=&\int_{L_1}^{L_2}[(H'(0))^2\psi_1(0,x)Z(\tau,x)-(H'(0))^2\psi_1(0,x)Z(\tau,x)]dt&\\[1mm]
&=&0
\end{array}$$
to conclude that
\begin{equation}
0\geq a\int_{0^+}^{\tau^-}(\langle \phi_1,Z\rangle-\langle W,\psi_1\rangle)dt+(\underline{\lambda}_1-\overline{\lambda}_1)\int_{0^+}^{\tau^-}\langle\psi_1,Z\rangle dt.
\label{c18}
\end{equation}

Finally, adding \eqref{c17} and \eqref{c18} yields
$$
(\underline{\lambda}_1-\overline{\lambda}_1)(\frac{1}{a}\int_{0^+}^{\tau^-}\langle\psi_1,Z\rangle dt+\frac{1}{b}\int_{0^+}^{\tau^-}\langle\phi_1,W\rangle dt)\leq0.
$$
The result $\underline{\lambda}_1\leq\overline{\lambda}_1$ is now proved since $\frac{1}{a}\int_{0^+}^{\tau^-}\langle\psi_1,Z\rangle dt+\frac{1}{b}\int_{0^+}^{\tau^-}\langle\phi_1,W\rangle dt\geq,\not\equiv 0$.

(ii) If $J_1=J_2$, then $\lambda^*((L_1,L_2),H'(0))$ is the principal eigenvalue to problem \eqref{c02} with $\mathcal{\tilde{L}}_1=\mathcal{\tilde{L}}_2(:=\mathcal{\tilde{L}})$. We first consider the following auxiliary problem
\begin{eqnarray}
\left\{
\begin{array}{lll}
-\phi^*_t-d_1\mathcal{\tilde{L}}[\phi^*]=a(t)\psi^*-[a(t)+m_1(t)]\phi^*+\lambda_2\phi^*,& t\in(0^+,\tau],\,\,x \in (L_1,L_2),  \\[1mm]
-\psi^*_t-d_2\mathcal{\tilde{L}}[\psi^*]=b(t)\phi^*-m_2(t)\psi^*+\lambda_2\psi^*, & t\in(0^+,\tau],\,\,x \in (L_1,L_2),  \\[1mm]
\phi^*(0^+,x)=\phi^*(0,x),& x\in (L_1,L_2),\\[1mm]
\psi^*(0^+,x)=\frac{1}{H'(0)}\psi^*(0,x),& x\in (L_1,L_2),\\[1mm]
\phi^*(0,x)=\phi^*(\tau,x),\psi^*(0,x)=\psi^*(\tau,x), & x\in [L_1,L_2].
\end{array} \right.
\label{c14}
\end{eqnarray}
Multiplying the first equation in \eqref{c02} by $\phi^*$ and the first equation in \eqref{c14} by $\phi$, then  integrating both sides of these two equations in $x$ over $[L_1,L_2]$ give that
$$
\langle\phi_t,\phi^*\rangle-d_1\langle\mathcal{\tilde{L}}[\phi],\phi^*\rangle=b(t)\langle \psi,\phi^*\rangle-\langle[a(t)+m_1(t)]\phi,\phi^*\rangle+\lambda^*\langle\phi,\phi^*\rangle,
$$
$$
\langle-\phi^*_t,\phi\rangle-d_1\langle\mathcal{\tilde{L}}[\phi^*],\phi\rangle=a(t)\langle \psi^*,\phi\rangle-\langle[a(t)+m_1(t)]\phi^*,\phi\rangle+\lambda_2\langle\phi^*,\phi\rangle.
$$
Abstracting the above two equations and then integrating in $t$ over $(0^+,\tau]$, asserts
\begin{eqnarray}
0=\int_{0^+}^{\tau}(b(t)\langle\psi,\phi^*\rangle-a(t)\langle \psi^*,\phi\rangle)dt+(\lambda^*-\lambda_2)\int_{0^+}^{\tau}\langle\phi,\phi^*\rangle dt
\label{c15}
\end{eqnarray}
in view of $\langle\mathcal{\tilde{L}}[\phi],\phi^*\rangle=\langle\mathcal{\tilde{L}}[\phi^*],\phi\rangle$ and the periodicity and impulsive conditions of $\phi$ and $\phi^*$.
Similarly, careful calculations yield
\begin{eqnarray}
0=\int_{0^+}^{\tau}(a(t)\langle\phi,\psi^*\rangle-b(t)\langle \phi^*,\psi\rangle)dt+(\lambda^*-\lambda_2)\int_{0^+}^{\tau}\langle\psi,\psi^*\rangle dt.
\label{c16}
\end{eqnarray}
By adding \eqref{c15} and \eqref{c16}, it is clear that $\lambda_2=\lambda^*$.

Next we take the problem \eqref{c14} with $(\lambda^*,\phi^*,\psi^*)$ and \eqref{c13} with $(\overline{\Lambda},\Phi,\Psi)$ into consideration. By the same procedure, we obtain
$$
0\leq\int_{0^+}^{\tau}(b(t)\langle\Psi,\phi^*\rangle-a(t)\langle \psi^*,\Phi\rangle)dt+(\overline{\Lambda}-\lambda^*)\int_{0^+}^{\tau}\langle\Phi,\phi^*\rangle dt
$$
$$
0\leq\int_{0^+}^{\tau}(a(t)\langle\Phi,\psi^*\rangle-b(t)\langle \phi^*,\Psi\rangle)dt+(\overline{\Lambda}-\lambda^*)\int_{0^+}^{\tau}\langle\Psi,\psi^*\rangle dt
$$
which concludes that $\lambda^*\leq \overline{\Lambda}$ since $\int_{0^+}^{\tau}(\langle\Phi,\phi^*\rangle+\langle\Psi,\psi^*\rangle)dt \geq,\not\equiv 0$. Similarly, using \eqref{c14} with $(\lambda^*,\phi^*,\psi^*)$ and \eqref{c133} with $(\underline{\Lambda},\Phi_1,\Psi_1)$ to obtain $\lambda^*\geq \underline{\Lambda}$, so $$\underline{\Lambda}\leq\lambda^*\leq \overline{\Lambda}.$$
\epf

\begin{cor}\label{Co3.1}
Suppose $J_1\neq J_2$ and all coefficients in \eqref{c02} are constant. For the generalized principal eigenvalue, we have $\underline{\lambda}\leq\overline{\lambda}$.
\end{cor}
\bpf
 It follows from the definition of $\underline{\lambda}$ that $\underline \lambda_1\geq \underline{\lambda}-\epsilon$ for sufficiently small $\epsilon$, then $\underline{\lambda}\geq \underline\lambda_1\geq \underline{\lambda}-\epsilon$.
Similarly, $\overline{\lambda}\leq\overline\lambda_1\leq \overline{\lambda}+\epsilon$ holds.
Therefore, we obtain
$$
\underline{\lambda}-\epsilon\leq\underline\lambda_1\leq\overline\lambda_1\leq \overline{\lambda}+\epsilon
$$
by Lemma \ref{le3.1} (i). And letting $\epsilon\rightarrow0$ yields $\underline{\lambda}\leq \overline{\lambda}$.
\epf

If $J_1\neq J_2$ and coefficients in \eqref{c02} are time-dependent, whether $\underline{\lambda}\leq \overline{\lambda}$ still holds or not is subject to further discussion.

\vspace{2mm}
The properties of the generalized principal eigenvalue and the principal eigenvalue involving harvesting pulse are investigated in the following.

\begin{lem} \label{le3.2}Suppose that $\mathbf{(J)}$ holds, we have the following statements.

$(i)$ $\underline{\lambda}((L_1,L_2),H'(0))=\underline{\lambda}((0,L_2-L_1),H'(0))$ and $\overline{\lambda}((L_1,L_2),H'(0))=\overline{\lambda}((0,L_2-L_1),H'(0))$.

$(ii)$ $\underline{\lambda}(\Omega,H'(0))$ and $\overline{\lambda}(\Omega,H'(0))$ are nonincreasing with respect to the domain $\Omega$ for any given $H'(0)$, that is,
$$\underline{\lambda}(\Omega_2,H'(0)) \leq \underline{\lambda}(\Omega_1,H'(0))$$
and
$$\overline{\lambda}(\Omega_2,H'(0)) \leq \overline{\lambda}(\Omega_1,H'(0))$$
hold for any $\Omega_1\subseteq \Omega_2\subseteq R^1$.

$(iii)$ Furthermore, if $J_1(x)=J_2(x)$, then for any given $L>0$, $\lambda^*((0,L),H'(0))$ is strictly decreasing in $H'(0)$, and $\lambda^*((0,L),H'(0))$ is strictly decreasing with respect to the length $L$ of the interval $(0,L)$ for any given $H'(0)$.
\end{lem}
\bpf
(i) Let $y=x-L_1$ with $x\in[L_1,L_2]$, then $\underline{\lambda}((L_1,L_2),H'(0))=\underline{\lambda}((0,L_2-L_1),H'(0))$ and $\overline{\lambda}((L_1,L_2),H'(0))=\overline{\lambda}((0,L_2-L_1),H'(0))$ can be guaranteed by the transformation $(t,x)\rightarrow(t,y)$.

(ii)Without loss of generality, we only prove that $\underline{\lambda}(\Omega_2,H'(0)) \leq \underline{\lambda}(\Omega_1,H'(0))$ for any $\Omega_1\subseteq \Omega_2\subseteq R^1$, and it suffices to prove $\underline{\lambda}((0,L_2-L_1),H'(0)) \leq \underline{\lambda}((0,L_4-L_3),H'(0))$ if $L_2-L_1>L_4-L_3$.

Let $(\underline{\lambda}((0,L_2-L_1),H'(0)),\phi(t,x),\psi(t,x))$ be a generalized eigen-pair of \eqref{c02} on $(0^+,\tau]\times[0,L_2-L_1]$, satisfying
\begin{eqnarray}
\left\{
\begin{array}{lll}
\phi_t-d_1\mathcal{\tilde{L}}_1[\phi]\geq b(t)\psi-[a(t)+m_1(t)]\phi+\underline{\lambda}\phi,& t\in(0^+,\tau],\,\,x \in (L_1,L_2),  \\[1mm]
\psi_t-d_2\mathcal{\tilde{L}}_2[\psi]\geq a(t)\phi-m_2(t)\psi+\underline{\lambda} \psi, & t\in(0^+,\tau],\,\,x \in (L_1,L_2),  \\[1mm]
\phi(0^+,x)=\phi(0,x),& x\in (L_1,L_2),\\[1mm]
\psi(0^+,x)=H'(0)\psi(0,x),& x\in (L_1,L_2),\\[1mm]
\phi(0,x)=\phi(\tau,x),\psi(0,x)=\psi(\tau,x), & x\in [L_1,L_2].
\end{array} \right.
\label{c08}
\end{eqnarray}
Since
$$\begin{array}{llllll}
&&-d_1[\int_{0}^{L_2-L_1}J_1(x-y)\phi(t,y)dy-\phi(t,x)]\\[1mm]
&=&-d_1[\int_{0}^{L_4-L_3}J_1(x-y)\phi(t,y)dy-\phi(t,x)]-d_1\int_{L_4-L_3}^{L_2-L_1}J_1(x-y)\phi(t,y)dy\\[1mm]
&\leq,\not\equiv&-d_1[\int_{0}^{L_4-L_3}J_1(x-y)\psi(t,y)dy-\phi(t,x)]\\[1mm]
\end{array}$$
for $t\in(0^+,\tau)$ and $x\in[0,L_4-L_3]$, we obtain from the first equation in \eqref{c08} that
\begin{equation}
\begin{array}{lll}
&&\phi_t-d_1[\int_{0}^{L_4-L_3}J_1(x-y)\phi(t,y)dy-\phi(t,x)]\\[1mm]
&\geq,\not\equiv& b(t)\psi-[a(t)+m_1(t)]\phi+\underline{\lambda}((0,L_2-L_1),H'(0))\phi.
\end{array}
\label{c188}
\end{equation}
Similarly, for $t\in(0^+,\tau)$ and $x\in[0,L_4-L_3]$, we have
\begin{equation}
\begin{array}{lll}
&&\psi_t-d_1[\int_{0}^{L_4-L_3}J_2(x-y)\psi(t,y)dy-\phi(t,x)]\\[1mm]
&\geq,\not\equiv& a(t)\phi-m_2(t)\psi+\underline{\lambda}((0,L_2-L_1),H'(0))\psi.
\label{c189}
\end{array}
\end{equation}

Furthermore,
$$\phi(0,x)=\phi(\tau,x),\,\psi(0,x)=\psi(\tau,x)$$
and
$$\phi(0^+,x)=\phi(0,x),\,\psi(0^+,x)=H'(0)\psi(0,x)$$
in $x\in[0,L_2-L_1]$. Finally, the definition of $\underline{\lambda}((0,L_4-L_3),H'(0))$ guarantees that $\underline{\lambda}((0,L_2-L_1),H'(0))\leq\underline{\lambda}((0,L_4-L_3),H'(0))$, So  $\underline{\lambda}(\Omega,H'(0))$ is nonincreasing in $\Omega$.

 (iii) Since $J_1(x)=J_2(x)$, the principal eigenvalue problem \eqref{c02} admits a unique solution $(\lambda^*((L_1,L_2),H'(0)),\phi,\psi)$ with $\phi,\psi>0$. In problem \eqref{c03} with constant coefficients replaced by time-periodic coefficients, differentiating both sides of the first two equations with respect to $H'(0)$ in $(0^+,\tau]$ yields
 \begin{eqnarray}
\left\{
\begin{array}{lll}
\alpha_t'=b(t)\beta'(t)-[a(t)+m_1(t)-d_1\lambda_0]\alpha'(t)+(\lambda^*)'\alpha(t)+\lambda^*\alpha'(t),& t\in(0^+,T],\\[1mm]
\beta_t'=a(t)\alpha'(t)-[m_2(t)-d_2\lambda_0]\beta'(t)+(\lambda^*)'\beta(t)+\lambda^*\beta'(t),& t\in(0^+,T],\\[1mm]
\alpha'(0)=\alpha'(\tau),\beta'(0)=\beta'(\tau),\\[1mm]
\alpha'(0^+)=\alpha'(0),\beta'(0^+)=\beta(0)+H'(0)\beta'(0),\\[1mm]
\end{array} \right.
\label{c09}
\end{eqnarray}
where $\alpha':=\frac{\partial\alpha(t)}{\partial H'(0)}$ and $\beta':=\frac{\partial\beta(t)}{\partial H'(0)}$ for brevity.

We now take the following auxiliary problem related to problem \eqref{c03} into consideration
\begin{eqnarray}
\left\{
\begin{array}{lll}
-\alpha^*_t=a(t)\beta^*(t)-[a(t)+m_1(t)-d_1\lambda_0]\alpha^*(t)+\lambda^*\alpha^*(t),& t\in(0^+,T],\\[1mm]
-\beta^*_t=b(t)\alpha^*(t)-[m_2(t)-d_2\lambda_0]\beta^*(t)+\lambda^*\beta^*(t),& t\in(0^+,T],\\[1mm]
\alpha^*(0)=\alpha^*(\tau),\beta^*(0)=\beta^*(\tau),\\[1mm]
\alpha^*(0^+)=\alpha^*(0),\beta^*(0^+)=\frac{1}{H'(0)}\beta^*(0).\\[1mm]
\end{array} \right.
\label{c10}
\end{eqnarray}
Multiplying the first equation in \eqref{c09} by $\alpha^*$ and the second equation by $\beta^*$, then integrating both sides in $=(0^+,T]$ yields yields
$$
-\int_{0^+}^{\tau}\alpha'\alpha^*_tdt=\int_{0^+}^{\tau}[b(t)\beta'\alpha^*-(a(t)+m_1(t)-d_1\lambda_0)\alpha'\alpha^*
+(\lambda^*)'\alpha\alpha^*+\lambda^*\alpha'\alpha^*]dt
$$
and
$$
-\frac{1}{H'(0)}\beta(0)\beta^*(0)-\int_{0^+}^{\tau}\beta'\beta^*_tdt=\int_{0^+}^{\tau}[a(t)\alpha'\beta^*-(m_2(t)-d_2\lambda_0)\beta'\beta^*
+(\lambda^*)'\beta\beta^*+\lambda^*\beta'\beta^*]dt.
$$
Recalling \eqref{c10}, it is easy to see that
\begin{eqnarray*}
\left\{
\begin{array}{lll}
\int_{0^+}^{\tau}[a(t)\alpha'\beta^*-b(t)\beta'\alpha^*]dt=\int_{0^+}^{\tau}(\lambda^*)'\alpha\alpha^*dt\\[1mm]
\int_{0^+}^{\tau}[b(t)\beta'\alpha^*-a(t)\alpha'\beta^*]dt-\frac{1}{H'(0)}\beta(0)\beta^*(0)=\int_{0^+}^{\tau}(\lambda^*)'\beta\beta^*dt\\[1mm]
\end{array} \right.
\end{eqnarray*}

Adding this two equations yields
$$
(\lambda^*)'((0,L),H'(0))
=\frac{-\frac{1}{H'(0)}\beta(0)\beta^*(0)}{\int_{0^+}^{\tau}(\beta\beta^*+\alpha\alpha^*)dt}<0,
$$
thus $\lambda^*((0,L),H'(0))$ is strictly decreasing in $H'(0)$ for any given $\Omega$.

Meanwhile, we can modify the proof of (ii) to derive the strictly decreasing of $\lambda^*((0,L),H'(0))$ in $L$. In fact, if $\underline{\lambda}=\overline{\lambda}:=\lambda^*$, then the inequalities in \eqref{c08} is replaced by equalities and $\underline{\lambda}((0,L_2-L_1),H'(0))$ is replaced by $\lambda^*(0,L_2-L_1),H'(0)$, respectively. Subsequently, it follows from inequalities in \eqref{c188} and \eqref{c189} with $\underline{\lambda}((0,L_2-L_1),H'(0))$ replaced by $\lambda^*((0,L_2-L_1),H'(0))$ and Lemma \ref{le3.1}(ii) that $\lambda^*((0,L_2-L_1),H'(0))<\lambda^*((0,L_4-L_3),H'(0))$, therefore $\lambda^*((0,L),H'(0))$ is strictly decreasing in $L$.
\epf


\section{Spreading and vanishing}
In order to investigate the criteria for species spread and vanish in impulsive and nonlocal model with moving boundaries, we first fix the domain from $[g(t),h(t)]$ to $[L_1,L_2]$.
\subsection{A corresponding fixed boundary problem}
The initial boundary problem with harvesting pulse and nonlocal diffusion in a fixed domain can be written by
\begin{eqnarray}
\left\{
\begin{array}{lll}
u_{1t}=d_1\mathcal{\tilde{L}}_1[u_1]+{b(t)}u_2-{a(t)}u_1-m_1(t)u_1\\[1mm]
\qquad-\alpha_1(t)u_1^2,  & t\in((n\tau)^+,(n+1)\tau],\,x \in [L_1,L_2],\\[1mm]
u_{2t}=d_2\mathcal{\tilde{L}}_2[u_2]+{a(t)}u_1-m_2(t)u_2-\alpha_2(t)u_2^2, & t\in((n\tau)^+,(n+1)\tau],\,x \in [L_1,L_2],\\[1mm]
u_1((n\tau)^+,x)=u_1(n\tau,x), & x \in (L_1,L_2), \\[1mm]
u_2((n\tau)^+,x)=H(u_2(n\tau,x)), & x \in (L_1,L_2), \\[1mm]
u_i(0,x)=u_{i,0}(x), & x \in [L_1,L_2],\,i=1,2,\\[1mm]
\end{array} \right.
\label{d01}
\end{eqnarray}
and its corresponding periodic problem is
\begin{eqnarray}
\left\{
\begin{array}{lll}
U_{1t}=d_1\mathcal{\tilde{L}}_1[U_1]+{b(t)}U_2-{a(t)}U_1-m_1(t)U_1-\alpha_1(t)U_1^2,  & t\in(0^+,\tau],\,x \in [L_1,L_2],\\[1mm]
U_{2t}=d_2\mathcal{\tilde{L}}_2[U_2]+{a(t)}U_1-m_2(t)U_2-\alpha_2(t)U_2^2, & t\in(0^+,\tau],\,x \in [L_1,L_2],\\[1mm]
U_1(0^+,x)=U_1(0,x), & x \in (L_1,L_2), \\[1mm]
U_2(0^+,x)=H(U_2(0,x)), & x \in (L_1,L_2), \\[1mm]
U_1(0,x)=U_1(\tau,x),\,U_1(0,x)=U_1(\tau,x), & x \in [L_1,L_2].\\[1mm]
\end{array} \right.
\label{d02}
\end{eqnarray}

\begin{defi}
$(\hat{u}_1(t,x),\hat{u}_2(t,x))\in [Q_n((0,+\infty)\times(L_1,L_2))]^2$ is called a lower solution to problem \eqref{d01} if
\begin{eqnarray*}
\left\{
\begin{array}{ll}
\hat u_{1t}\leq d_1\mathcal{\tilde{L}}_1[\hat u_1]+{b(t)}\hat u_2-{a(t)}\hat u_1-m_1(t)\hat u_1-\alpha_1(t)\hat u_1^2,  & t\in((n\tau)^+,(n+1)\tau],\,x \in [L_1,L_2],\\[1mm]
\hat u_{2t}\leq d_2\mathcal{\tilde{L}}_2[\hat u_2]+{a(t)}\hat u_1-m_2(t)\hat u_2-\alpha_2(t)\hat u_2^2, & t\in((n\tau)^+,(n+1)\tau],\,x \in [L_1,L_2],\\[1mm]
\hat u_1((n\tau)^+,x)\leq u_1(\hat n\tau,x), & x \in (L_1,L_2), \\[1mm]
\hat u_2((n\tau)^+,x)\leq H(\hat u_2(n\tau,x)), & x \in (L_1,L_2), \\[1mm]
\hat u_i(0,x)\leq \hat u_{i,0}(x), & t\in(0,+\infty),\,x \in [L_1,L_2],\,i=1,2\\[1mm]
\end{array} \right.
\end{eqnarray*}
with
$$
Q_n([0,+\infty)\times[L_1,L_2])\triangleq\{(u_1(t,x),u_2(t,x))\in C((n\tau,(n+1)\tau]\times[L_1,L_2])\}.
$$
Similarly, $(\tilde{u}_1(t,x),\tilde{u}_2(t,x))$ is called an upper solution to problem \eqref{d01} if all inequalities in the above are reversed. If also $(\hat{u}_1(t,x),\hat{u}_2(t,x))\leq(\tilde{u}_1(t,x),\tilde{u}_2(t,x))$, then $(\tilde u_1,\tilde u_2)$ and $(\hat u_1,\hat u_2)$ are called an ordered upper and lower solution. By this manner, the lower and upper solution $(\hat{U}_1^*(t,x),\hat{U}_2^*(t,x))$ and $(\tilde{U}_1^*(t,x),\tilde{U}_2^*(t,x))$ to  periodic problem \eqref{d02} can also be defined with initial value condition $\hat u_i(0,x)\leq\hat u_{i,0}(x)(\tilde u_i(0,x)\geq\tilde u_{i,0}(x))$ replaced by periodic condition $\hat U_i(0,x)\leq\hat U_i(\tau,x)(\tilde U_i(0,x)\geq\tilde U_i(\tau,x))$, respectively.
\end{defi}

\begin{thm}\label{th4.1}
Assume that $\mathbf{(J)}$ holds. Let $\overline{\lambda}((L_1,L_2),H'(0))$ and $\underline{\lambda}((L_1,L_2),H'(0))$ be the generalized principal eigenvalue of problem \eqref{c02} and  $(u_1(t,x),u_2(t,x))$ be the unique solution to problem \eqref{d01}. Then we have the following three statements.\\
$(i)$ The periodic problem \eqref{d02} admits a unique positive solution if $\overline{\lambda}((L_1,L_2),H'(0))<0$. Moreover, if $\underline{\lambda}((L_1,L_2),H'(0))\geq0$, $(0,0)$ is the only nonnegative periodic solution to problem \eqref{d02} provided that $J_1\neq J_2$ and all coefficients are constant.\\
$(ii)$ If $\underline{\lambda}((L_1,L_2),H'(0))\geq0$, then $\lim\limits_{t\to+\infty}{(u_1(t,x),u_2(t,x))}=(0,0)$ uniformly for $x\in[L_1,L_2]$.\\
$(iii)$ If $\overline{\lambda}((L_1,L_2),H'(0))<0$, then $\lim\limits_{m\to+\infty}{(u_1,u_2)}(t+m\tau,x)=(U_1^*,U_2^*)(t,x)$ uniformly for $(t,x)\in[0,\tau]\times[L_1,L_2]$, where $(U_1^*(t,x),U_2^*(t,x))$ is the unique positive solution to periodic problem \eqref{d02}.
\end{thm}
\bpf (i) Since $\overline{\lambda}((L_1,L_2),H'(0))<0$, one easily checks that $(A,A)$ and $(\underline{U}_1,\underline{U}_2)$ are the upper and lower solutions to problem \eqref{d02}, where $(A,A)$ is defined in Lemma \ref{m2} and $(\underline{U}_1,\underline{U}_2)$ satisfies
\begin{equation*}
\underline{U}_1(t,x)=\varepsilon e^{(\lambda+\upsilon)\tau}\phi(t,x),\ t\in[0,+\infty),\,x\in[L_1,L_2]
\end{equation*}
and
\begin{eqnarray*}
\underline{U}_2(t,x)=
\begin{cases}
\displaystyle \varepsilon{\psi}(n\tau,x),&\,\,x\in[L_1,L_2],\\[3mm]
\displaystyle \varepsilon e^{(\lambda+\upsilon)\tau}\psi((n\tau)^+,x), &\,\,x\in[L_1,L_2],\\[3mm]
\displaystyle \varepsilon e^{(\lambda+\upsilon)\tau}e^{(-\lambda-\upsilon)(t-n\tau)}\psi(t,x),&\,\,t\in((n\tau)^+,(n+1)\tau],\,\,x\in[L_1,L_2]
\end{cases}
\end{eqnarray*}
for suitable $\upsilon$ \cite[Theorem 4.3]{XLS}, where $(\lambda, \phi, \psi)$ satisfies \eqref{c01}. Then the iteration sequence $\{\overline U_1^{(m)},\overline U_2^{(m)}\}$  with $(\overline U_1^{(0)},\overline U_2^{(0)})=(A,A)$  can be constructed, which satisfy
\begin{eqnarray}
\left\{
\begin{array}{lll}
\overline U_{1t}^{(m)}-d_1\mathcal{\tilde{L}}_1[\overline U_1^{(m)}]+K_1\overline U_1^{(m)}={b(t)}\overline U_2^{(m-1)}-{a(t)}\overline U_1^{(m-1)}\\[1mm]
\qquad \qquad\quad-m_1(t)\overline U_1^{(m-1)}-\alpha_1(t)(\overline U_1^{(m-1)})^2+K_1\overline U_1^{(m-1)},  & t\in(0^+,\tau],\,x \in [L_1,L_2],\\[1mm]
\overline U_{2t}^{(m)}-d_2\mathcal{\tilde{L}}_2[\overline U_2^{(m)}]+K_2\overline U_2^{(m)}={a(t)}\overline U_1^{(m-1)}-m_2(t)\overline U_2^{(n-1)}\\[1mm]
\qquad\qquad \quad-\alpha_2(t)(\overline U_2^{(m-1)})^2+K_2\overline U_2^{(m-1)}, & t\in(0^+,\tau],\,x \in [L_1,L_2],\\[1mm]
\overline U_1^{(m)}(0^+,x)=\overline U_1^{(m-1)}(\tau,x), & x \in (L_1,L_2), \\[1mm]
\overline U_2^{(m)}(0^+,x)=H(\overline U_2^{(m-1)}(\tau,x)), & x \in (L_1,L_2), \\[1mm]
\overline U_i^{(m)}(0,x)=\overline U_i^{(m-1)}(\tau,x), & t\in(0,+\infty),\,x \in [L_1,L_2],\\[1mm]
\end{array} \right.
\label{ff1}
\end{eqnarray}
where $K_1$ and $K_2$ are large enough. Similarly, we can define iteration sequence $\{\underline U_1^{(m)},\underline U_2^{(m)}\}$ with $(\underline U_1^{(0)},\underline U_2^{(0)})=(\underline U_1,\underline U_2)$.

We claim that sequences $\{\overline U_1^{(m)},\overline U_2^{(m)}\}$ and $\{\underline U_1^{(m)},\underline U_2^{(m)}\}$ are monotonically decreasing and increasing, respectively, and they converge to positive function pairs $(\overline U_1^*(t,x),\overline U_2^*(t,x))$ and $(\underline U_1^*(t,x),\underline U_2^*(t,x))$. Furthermore, $(\overline U_1^*(t,x),\overline U_2^*(t,x))$ and $(\underline U_1^*(t,x),\underline U_2^*(t,x))$ are the maximal and minimal positive periodic solution to problem \eqref{d02}. Moreover,
\begin{equation}
(\underline{U}_1,\underline{U}_2)\leq(\overline U_1^*(t,x),\overline U_2^*(t,x)),(\underline U_1^*(t,x),\underline U_2^*(t,x))\leq (A,A).
\label{ff}
\end{equation}

The uniqueness of positive periodic solution $(U_1^*,U_2^*)$ will be shown in the following. To proceed, let $(U_1^*,U_2^*)$ and $(\tilde{U}_1,\tilde{U}_2)$ be two different solutions, and denote
$$
S=\{s\in[0,1],s\tilde{U}_1\leq{U_1^*},\,s\tilde{U}_2\leq{U_2^*},\,t\in[0,T],\,x\in[L_1,L_2]\}.
$$
We claim that $1\in S$. Otherwise, we assume that $s_0=\sup S<1$. Define $(W_1,W_2):=(U_1^*-s_0\tilde{U}_1,U_2^*-s_0\tilde{U}_2)$, satisfying
$(W_1,W_2)\geq(0,0)$ in $[0,\tau]\times [L_1,L_2]$. Without loss of generality, we assume that $W_1(t_0,x_0)=0$ for some $(t_0,x_0)\in [0,\tau]\times [L_1,L_2]$. Direct calculations yield
$$\begin{array}{llllll}
&&W_{1t}-d_1\mathcal{\tilde{L}}_1[W_1]\\[1mm]
&=&{b(t)}(U_2^*-s_0\tilde{U}_2)-({a(t)}+m_1(t))(U_1^*-s_0\tilde{U}_1)-\alpha_1(t)[(U_1^*)^2-s_0\tilde{U}_1^2]\\[1mm]
&>&-({a(t)}+m_1(t))(U_1^*-s_0\tilde{U}_1)-\alpha_1(t)(U_1^*-s_0\tilde{U}_1)(U_1^*+s_0\tilde{U}_1)\\[1mm]
&=&-K_1(t)W_1\\[1mm]
\end{array}$$
with $K_1(t):={a(t)}+m_1(t)+\alpha_1(t)(U_1^*+s_0\tilde{U}_1)$ and
$$\begin{array}{llllll}
&&W_{2t}-d_2\mathcal{\tilde{L}}_2[W_2]\\[1mm]
&=&{a(t)}(U_1^*-s_0\tilde{U}_1)-m_2(t)(U_2^*-s_0\tilde{U}_2)-\alpha_2(t)[(U_2^*)^2-s_0\tilde{U}_2^2]\\[1mm]
&>&-m_2(t)(U_2^*-s_0\tilde{U}_2)-\alpha_2(t)(U_2^*-s_0\tilde{U}_2)(U_2^*+s_0\tilde{U}_2)\\[1mm]
&=&-K_2(t)W_2\\[1mm]
\end{array}$$
with $K_2(t):=m_2(t)+\alpha_2(t)(U_2^*+s_0\tilde{U}_2).$ So $$W_{1t}(t_0,x_0)>d_1\mathcal{\tilde{L}}_1[W_1]=d_1\int_{L_1}^{L_2}J_1(x_0-y)W_1(t_0,y)dy\geq0.$$ On the other hand, $W_{1t}(t_0,x_0)\leq0$ since $W_1(t_0,x_0)=0$ for some $(t_0,x_0)\in [0,\tau]\times [L_1,L_2]$, which leads contradiction, so the uniqueness of solution to problem \eqref{d02} is now obtained in $[0,\tau]\times[L_1,L_2]$. Since $(A,A)$ and $(\underline{U}_1,\underline{U}_2)$ can be sufficiently large and small, respectively, it then follows from \eqref{ff} that the solution $(U_1^*(t,x),U_2^*(t,x))$ to periodic problem \eqref{d02} is unique in the whole space.

We finally claim that problem \eqref{d02} has no positive solution if $\underline{\lambda}((L_1,L_2),H'(0))\geq0$. Suppose by contradiction that $(\overline{U}_1,\overline{U}_2)$ is a positive steady state of \eqref{d02}, which meets
\begin{eqnarray*}
\left\{
\begin{array}{lll}
\overline U_{1t}-d_1\mathcal{\tilde{L}}_1[\overline U_1]={b(t)}\overline U_2-{a(t)}\overline U_1-m_1(t)\overline U_1-\alpha_1(t)\overline U_1^2\\[1mm]
\qquad \qquad\qquad\quad<{b(t)}\overline U_2-{a(t)}\overline U_1-m_1(t)\overline U_1,  & t\in(0^+,\tau],\,x \in [L_1,L_2],\\[1mm]
\overline U_{2t}-d_2\mathcal{\tilde{L}}_2[\overline U_2]={a(t)}\overline U_1-m_2(t)\overline U_2-\alpha_2(t)\overline U_2^2\\[1mm]
\qquad \qquad\qquad\quad<{a(t)}\overline U_1-m_2(t)\overline U_2, & t\in(0^+,\tau],\,x \in [L_1,L_2],\\[1mm]
\overline U_1(0^+,x)=\overline U_1(0,x), & x \in (L_1,L_2), \\[1mm]
\overline U_2(0^+,x)=H(\overline U_2(0,x)), & x \in (L_1,L_2), \\[1mm]
\overline U_1(0,x)=\overline U_1(\tau,x),\,\overline U_1(0,x)=\overline U_1(\tau,x), &x \in [L_1,L_2].
\end{array} \right.
\end{eqnarray*}
It is deduced by the definition of generalized principal eigenvalue that $\overline{\lambda}((L_1,L_2),H'(0))\leq0$. We have the following two cases.

\textbf{Case 1} If $\overline{\lambda}=\underline{\lambda}$, then the principal eigenvalue $\lambda^*$ is well defined and $\lambda^*<0$ by the comparison principle in Lemma \ref{le3.1} (ii) with $\overline{\Lambda}=0$, which leads contradiction to $\lambda^*=\underline{\lambda}\geq0$.

\textbf{Case 2} If $\overline{\lambda}\neq\underline{\lambda}$, recalling Corollary \ref{Co3.1}, we have $\underline{\lambda}<\overline{\lambda}\leq0$, which also contradicts to the condition $\underline{\lambda}\geq0$.\\
Thus periodic problem \eqref{d02} has the unique nonnegative solution $(0,0)$ if $\underline{\lambda}((L_1,L_2),H'(0))\geq0$.

(ii) It is easy to verify that $(Me^{-\underline{\lambda}t}\phi(t,x),Me^{-\underline{\lambda}t}\psi(t,x))$ is an upper solution to problem \eqref{d01}, in which $(\phi,\psi)$ satisfy  problem \eqref{c01} and $M$ is a positive constant which is sufficient large such that $M\phi(0,x)\geq u_{1,0}(x)$ and $M\psi(0,x)\geq u_{2,0}(x)$ in $x\in[L_1,L_2]$. Since $(Me^{-\underline{\lambda}t}\phi(t,x),Me^{-\underline{\lambda}t}\psi(t,x))\rightarrow(0,0)$ as $t\rightarrow\infty$, it follows from the comparison principle that $\lim\limits_{t\to+\infty}{(u_1(t,x),u_2(t,x))}=(0,0)$ uniformly for $x\in[L_1,L_2]$.

(iii) We prove it by induction. Without loss of generality, we suppose that $u_2(0,x)>0$ in $x\in (L_1,L_2)$, otherwise, $u(\tau,x)$ can be regarded as the new initial value. In fact, for $m=0$, since $\epsilon$ is sufficiently small and $A$ is large enough, we obtain
$\underline{U}_2^{(0)}(0,x)\leq{u_2(0,x)}\leq{\overline{U}_2^{(0)}(0,x)}$ for $x\in [L_1,L_2]$, so $H(\underline{U}_2^{(0)}(0,x))\leq{H(u_2(0,x))}\leq{H(\overline{U}_2^{(0)}(0,x)})$ for $x\in [L_1,L_2]$, which means $\underline{U}_2^{(0)}(0^+,x)\leq{u_2(0^+,x)}\leq{\overline{U}_2^{(0)}(0^+,x)}$ for $x\in [L_1,L_2]$. Comparison principle yields $\underline{U}_2^{(0)}(t,x)\leq{u_2(t,x)}\leq{\overline{U}_2^{(0)}(t,x)}$ for $t\in[0,T]$ and $x\in [L_1,L_2]$. By the same procedure for $t\in[n\tau,(n+1)\tau](n\geq1)$, we finally obtain
$$
\underline{U}_2^{(0)}(t,x)\leq{u_2(t,x)}\leq{\overline{U}_2^{(0)}(t,x)},\, t\geq0,\,x\in [L_1,L_2].
$$
Similarly, for $m=1$, using the iteration sequence in \eqref{ff1} to obtain
$$
\underline{U}_2^{(1)}(t,x)\leq{u_2(t+\tau,x)}\leq{\overline{U}_2^{(1)}(t,x)},\, t\geq0,\,x\in [L_1,L_2],
$$
By induction of $m$, we finally get
$$
(\underline U_1^{(m)},\underline U_2^{(m)})(t,x)\leq (u_1,u_2)(t+m\tau,x)\leq(\overline U_1^{(m)},\overline U_2^{(m)})(t,x),\,t\geq0,x\in[L_1,L_2],
$$
which together with the uniqueness of periodic solution to \eqref{d02} in (i) by $\overline{\lambda}((L_1,L_2),H'(0))<0$, yields
 $\lim\limits_{m\to+\infty}{(u_1,u_2)}(t+m\tau,x)=(U_1^*,U_2^*)(t,x)$ uniformly for $x\in\overline[L_1,L_2]$.
\epf

\subsection{Main results for moving boundary problem}
This section is devoted to researching the criteria and sufficient conditions for spreading-vanishing of individuals, in which the negative effect of harvesting rate on the persistence of invasive species can be observed, and difficulties induced by nonlocal diffusion in theoretical
analysis are solved.

In view of $g'(t)<0$ and $h'(t)>0$, we declare
$$
\lim\limits_{t\to+\infty}g(t)=g_\infty\,\,\text{and}\,\,\lim\limits_{t\to+\infty}h(t)=h_\infty,
$$
with some $g_\infty\in [-\infty,-h_0)$ and $h_\infty\in (h_0,+\infty]$.

Since $-g(t)$ and $h(t)$ are strictly increasing with respect to $t$, it then follows from Lemma 3.2 (ii) that $\overline{\lambda}((g(t),h(t)),H'(0))$ and $\underline{\lambda}((g(t),h(t)),H'(0)))$ are nonincreasing in $t$, and $\lambda^*((g(t),h(t)),H'(0))$ is strictly decreasing in $t$ if the principal eigenvalue exists. We then define
$$
\overline{\lambda}((g_\infty,h_\infty),H'(0)):=\lim\limits_{t\to\infty}\overline{\lambda}((g(t),h(t)),H'(0))
$$
and
$$
\overline{\lambda}((-\infty,+\infty),H'(0)):=\lim\limits_{-L_1,L_2\to\infty}\overline{\lambda}((L_1,L_2),H'(0)),
$$
which all are related to pulse function $H'(0)$. Similarly, we can also define $\underline{\lambda}((g_\infty,h_\infty),H'(0))$, $\underline{\lambda}((-\infty,+\infty),H'(0))$, ${\lambda}^*((g_\infty,h_\infty),H'(0))$ and
${\lambda}^*((-\infty,+\infty),H'(0))$.

In the following, we aim to characterize the spreading or vanishing of juveniles and adults with harvesting pulse in nonlocal diffusion and free boundaries. To begin with, we first introduce the comparison principle, and its proof is standard, see \cite[Theorem 3.2]{CD} without pulse.
\begin{lem}\label{le4.2}
{\rm (Comparison principle)} Assume that $\mathbf{(J)}$ holds and $(u_1,u_2;(g,h))$ is the unique solution to problem \eqref{a01}. Let
$(\bar{u}_1(t,x),\bar{u}_2(t,x))\in [C((n\tau,(n+1)\tau]\times[\bar{g}(t),\bar{h}(t)])]^2$ for $n=0,1,2,\dots$, $(\bar{g}(t),\bar{h}(t))\in [C((0,+\infty))\cap C^1((n\tau,(n+1)\tau])]^2$, satisfying
\begin{eqnarray*}
\left\{
\begin{array}{lll}
\bar{u}_{1t}\geq d_1\mathcal{L}_1[\bar{u}_1]+{b(t)}\bar{u}_2-{a(t)}\bar{u}_1-m_1(t)\bar{u}_1\\[1mm]
\qquad-\alpha_1(t)\bar{u}_1^2,
& t\in((n\tau)^+,(n+1)\tau],\,x \in (\bar g(t),\bar h(t)),\\[1mm]
\bar{u}_{2t}\geq d_2\mathcal{L}_2[\bar{u}_2]+{a(t)}\bar{u}_1-m_2(t)\bar{u}_2-\alpha_2(t)\bar{u}_2^2,
& t\in((n\tau)^+,(n+1)\tau],\,x \in (\bar g(t),\bar h(t)),\\[1mm]
\bar{u}_1((n\tau)^+,x)\geq\bar{u}_1(n\tau,x), & x \in (\bar g(n\tau),\bar h(n\tau)), \\[1mm]
\bar{u}_2((n\tau)^+,x)\geq H(\bar{u}_2(n\tau,x)), & x \in (g(n\tau),h(n\tau)), \\[1mm]
\bar{u}_1(t,\bar g(t)),\bar{u}_1(t,\bar h(t)),\bar{u}_2(t,\bar g(t)),\bar{u}_2(t,\bar h(t))\geq0, & t\in(0,+\infty), \\[1mm]
\bar{h}'(t)\geq \sum\limits_{i=1}^{2}\mu_i\int_{\bar{g}(t)}^{\bar{h}(t)}\int_{\bar{h}(t)}^{+\infty}J_i(x-y)\bar{u}_i(t,x)dydx, &  t\in(n\tau,(n+1)\tau], \\[1mm]
\bar{g}'(t)\leq-\sum\limits_{i=1}^{2}\mu_i\int_{\bar{g}(t)}^{\bar{h}(t)}\int_{-\infty}^{\bar{g}(t)}J_i(x-y)\bar{u}_i(t,x)dydx, &  t\in(n\tau,(n+1)\tau], \\[1mm]
-\bar{g}(0)\leq h_0<h_0\leq\bar{h}(0), \\[1mm]
\bar{u}_i(0,x)=u_{i,0}(x), & x\in[-h_0,h_0],\,i=1,2,\\[1mm]
\end{array} \right.
\end{eqnarray*}
then $[g(t),h(t)]\subset[\bar{g}(t),\bar{h}(t)]$ for $t\in(0,+\infty)$ and
$$(u_1(t,x),u_2(t,x))\leq(\bar{u}_1(t,x),\bar{u}_2(t,x))\,\,\text{for}\,\, t\in(0,+\infty)\,\,\text{and}\,\,x\in[g(t),h(t)].$$
The triplet $(\bar{u}_1(t,x),\bar{u}_2(t,x);(\bar{g}(t),\bar{h}(t)))$ is called an upper solution to problem \eqref{a01}. Also, if all inequalities in the above are reversed and $(\bar{u}_1,\bar{u}_2;(\bar{g},\bar{h}))$ is replaced by $(\underline{u}_1,\underline{u}_2;(\underline g,\underline h))$, then $[g(t),h(t)]\supset [\underline{g}(t),\underline{h}(t)]$ and
$$(u_1(t,x),u_2(t,x))\geq(\underline{u}_1(t,x),\underline{u}_2(t,x))\,\,\text{for}\,\, t\in(0,+\infty)\,\,\text{and}\,\,x\in[\underline g(t),\underline h(t)].$$
\end{lem}

\begin{lem}\label{le4.3}
Assume that $J_1(x)=J_2(x)$. If $h_\infty-g_\infty<\infty$, then $\lambda^*((g_\infty,h_\infty),H'(0))\geq0$ and
$$\lim\limits_{t\to+\infty}{||(u_1(t,\cdot),u_2(t,\cdot))||_{C[g(t),h(t)]}}=(0,0)$$
uniformly for $x\in[g(t),h(t)]$.
\end{lem}
\bpf
Suppose by contradiction that $\lambda^*((g_\infty,h_\infty),H'(0))<0$. Then there exists a sufficiently small $\varepsilon_1$ such that
$$\lambda^*((g_\infty+\varepsilon,h_\infty-\varepsilon),H'(0))<0\,\,\text{for\,all}\,0<\varepsilon<\varepsilon_1.$$
 Also, by continuous dependence of $\varepsilon$, for such $\varepsilon$ we can find a large integer $N_\varepsilon$ such that
$$
g_\infty<g(t)<g_\infty+\varepsilon,\ h_\infty-\varepsilon<h(t)<h_\infty, \forall\,\, t\geq N_\varepsilon\tau.
$$
Considering the following initial value problem for $n\geq N_\varepsilon$
\begin{eqnarray*}
\left\{
\begin{array}{lll}
\underline{u}_{1t}=d_1(\int_{g_\infty+\varepsilon}^{h_\infty-\varepsilon}J_1(x-y)\underline{u}_1(t,y)dy-\underline{u}_1(t,x))\\[1mm]
\qquad+{b(t)}\underline{u}_2-(a(t)+m_1(t))\underline{u}_1
-\alpha_1(t)\underline{u}_1^2,
& (t,x)\in\Omega_{\tau,\varepsilon},\\[1mm]
\underline{u}_{2t}=d_2(\int_{g_\infty+\varepsilon}^{h_\infty-\varepsilon}J_2(x-y)\underline{u}_2(t,y)dy-\underline{u}_2(t,x))\\[1mm]
\qquad+{a(t)}\underline{u}_1-m_2(t)\underline{u}_2-\alpha_2(t)\underline{u}_2^2,
& (t,x)\in\Omega_{\tau,\varepsilon},\\[1mm]
\underline{u}_1((n\tau)^+,x)=\underline{u}_1(n\tau,x), & x \in (g_\infty+\varepsilon,h_\infty-\varepsilon), \\[1mm]
\underline{u}_2((n\tau)^+,x)= H(\underline{u}_2(n\tau,x)), & x \in (g_\infty+\varepsilon,h_\infty-\varepsilon), \\[1mm]
\underline{u}_i(N_{\varepsilon}\tau,x)=u_i(N_{\varepsilon}\tau,x), & x\in(g_\infty+\varepsilon,h_\infty-\varepsilon),\,i=1,2,\\[1mm]
\end{array} \right.
\end{eqnarray*}
where $\Omega_{\tau,\varepsilon}:=\{(t,x):(n\tau)^+<t\leq(n+1)\tau],\,g_\infty+\varepsilon<x<h_\infty-\varepsilon\}$.

It follows from the comparison principle in Lemma \ref{le4.2} that $(u_1,u_2)(t,x)\geq(\underline{u}_1,\underline{u}_2)(t,x)$ for $t\geq N_{\varepsilon}\tau$ and $x\in[g_\infty+\varepsilon,h_\infty-\varepsilon]$. Recalling $\lambda^*((g_\infty+\varepsilon,h_\infty-\varepsilon),H'(0))<0$ and Theorem \ref{th4.1}(iii), one easily checks that
$$\lim\limits_{n\to+\infty}{(\underline u_1(t+n\tau,x),\underline u_2(t+n\tau,x))}=(U_{1}^\varepsilon(t,x),U_{2}^\varepsilon(t,x)),\,t\in[0,\tau],\, x\in[g_\infty,h_\infty],$$
where$(U_1^\varepsilon,U_2^\varepsilon)$ is the unique positive steady state with initial value condition $\underline{u}_i(N_{\varepsilon}\tau,x)=u_i(N_{\varepsilon}\tau,x)$ replaced by periodic condition $\underline{u}_i(0,x)=\underline{u}_i(\tau,x)$ in $x\in(g_\infty+\varepsilon,h_\infty-\varepsilon)$.

Consequently, there exists a $N_{1\varepsilon}(\geq N_\varepsilon)$ such that
$$
(u_1,u_2)(t+n\tau,x)\geq(\underline{u}_1,\underline{u}_2)(t+n\tau,x)\geq\frac{1}{2}(U_1^\varepsilon,U_2^\varepsilon)(t,x)>0,\,\forall t\geq N_{1\varepsilon}\tau, x\in[g_\infty+\varepsilon,h_\infty-\varepsilon].
$$
Since $J_1(0)=J_2(0)>0$, $J_1(x)=J_2(x)>0$ for $x\in[-3\varepsilon,3\varepsilon]$ with such $\varepsilon\in(0,h_0)$, we then define
$$
\gamma=\min_{[-3\varepsilon,3\varepsilon]}J_1(x)\,\,\,\text{and}\,\,\, \,\,\delta=\min_{[0,\tau]\times[g_\infty+\varepsilon,h_\infty-\varepsilon]}\{U_1^\varepsilon(t,x),U_2^\varepsilon(t,x)\},
$$
then $\gamma$,$\delta>0$. Notice the fact $[h_\infty-2\varepsilon,h_\infty-\varepsilon]\subset [g_\infty+\varepsilon,h_\infty-\varepsilon]$, we obtain
$$\begin{array}{llllll}
&&h'(t+n\tau)=\sum\limits_{i=1}^{2}\mu_i\int_{g(t+n\tau)}^{h(t+n\tau)}\int_{h(t+n\tau)}^{+\infty}J_i(x-y)u_i(t+n\tau,x)dydx\\[1mm]
&\geq&\sum\limits_{i=1}^{2}\mu_i\int_{h(t+n\tau)-2\varepsilon}^{h(t+n\tau)}\int_{h(t+n\tau)}^{h(t+n\tau)
+\varepsilon}J_i(x-y)u_i(t+n\tau,x)dydx\\[1mm]
&\geq&\gamma\varepsilon\sum\limits_{i=1}^{2}\mu_i\int_{h(t+n\tau)-2\varepsilon}^{h(t+n\tau)}u_i(t+n\tau,x)dx\\[1mm]
&\geq&\gamma\varepsilon\sum\limits_{i=1}^{2}\mu_i\int_{h_\infty-2\varepsilon}^{h_\infty-\varepsilon}u_i(t+n\tau,x)dx,\\[1mm]
&\geq&\frac{1}{2}\gamma\delta\varepsilon^2\sum\limits_{i=1}^{2}\mu_i>0,\,\text{for}\,t\geq N_{1\varepsilon}\tau.\\[1mm]
\end{array}$$
This contradicts to $h_\infty<\infty$. Therefore $\lambda^*((g_\infty,h_\infty),H'(0))\geq0$.

The proof ends with $(u_1,u_2)(t,x)\rightarrow(0,0)$. In fact, we can construct an upper solution to problem \eqref{a01}, satisfying
\begin{eqnarray*}
\left\{
\begin{array}{lll}
\bar{u}_{1t}=d_1(\int_{g_\infty}^{h_\infty}J_1(x-y)\bar{u}_1(t,y)dy-\bar{u}_1(t,x))+{b(t)}\bar{u}_2\\[1mm]
\qquad-{a(t)}\bar{u}_1-m_1(t)\bar{u}_1-\alpha_1(t)\bar{u}_1^2,
& (t,x)\in\Omega_{\tau,\infty},\\[1mm]
\bar{u}_{2t}=d_2(\int_{g_\infty}^{h_\infty}J_2(x-y)\bar{u}_2(t,y)dy-\bar{u}_2(t,x))\\[1mm] \qquad+{a(t)}\bar{u}_1-m_2(t)\bar{u}_2-\alpha_2(t)\bar{u}_2^2,
& (t,x)\in\Omega_{\tau,\infty},\\[1mm]
\bar{u}_1((n\tau)^+,x)=\bar{u}_1(n\tau,x), & x \in (g_\infty,h_\infty), \\[1mm]
\bar{u}_2((n\tau)^+,x)= H(\bar{u}_2(n\tau,x)), & x \in (g_\infty,h_\infty), \\[1mm]
\bar{u}_i(0,x)=\widetilde{u}_{i,0}(x), & x\in[-h_0,h_0],\,i=1,2
\end{array} \right.
\end{eqnarray*}
for $\,n=0,1,2,\dots$, where $\Omega_{\tau,\infty}:=\{(t,x):(n\tau)^+<t\leq(n+1)\tau],\,g_\infty<x<h_\infty\}$, $\widetilde{u}_{i,0}(x)=u_{i,0}(x)$ in $[-h_0,h_0]$ and $\widetilde{u}_{i,0}(x)=0$ in $[g_\infty,-h_0)\bigcup (h_0,h_\infty]$. Then $(u_1,u_2)(t,x)\leq(\bar u_1,\bar u_2)(t,x)$ in $[0,\infty)\times [g(t),h(t)]$ by the comparison principle. Since $\lambda^*((g_\infty,h_\infty),H'(0))\geq0$, it follows from Theorem \ref{th4.1}(ii) that
$(\bar u_1,\bar u_2)(t,x)\rightarrow(0,0)$ in $C([g_\infty,h_\infty])$, which ends the proof.
\epf

\begin{lem}\label{le4.4}
Suppose $J_1(x)\neq J_2(x)$ and all coefficients are constant. If $\underline{\lambda}(\infty,H'(0))\geq0$, then vanishing occurs, that is, the solution $(u_1,u_2)$ to problem \eqref{a01} satisfies
$$\lim\limits_{t\to+\infty}{(u_1,u_2)(t,x)}=(0,0)$$
uniformly in $x\in[g(t),h(t)]$.
\end{lem}
\bpf
Let $(\bar u_1,\bar u_2)(t)$  be the solution to
\begin{eqnarray*}
\left\{
\begin{array}{lll}
{(\bar u_1)}_t=b(t)\bar u_2-(a(t)+m_1(t)) \bar u_1-\alpha_1(t)(\bar u_1)^2,  & t\in((n\tau)^+,(n+1)\tau],\\[1mm]
{(\bar u_2)}_t=a(t)\bar u_1-m_2(t) \bar u_2-\alpha_2(t) (\bar u_2)^2, & t\in((n\tau)^+,(n+1)\tau],\\[1mm]
\bar u_1((n\tau)^+)=\bar u_1(n\tau),\bar u_2((n\tau)^+)=H(\bar u_2(n\tau)),\\[1mm]
\bar u_1(0)=\bar u_2(0)=A,\\[1mm]
\end{array} \right.
\end{eqnarray*}
where $A$ is defined in Lemma \ref{m2}. Recalling $\int_{-\infty}^{+\infty}J_i(x)dx=1$, it is clear that
$$
d_i\mathcal{L}_i[\bar u_i]=d_i\int_{g(t)}^{h(t)}J_i(x-y)\bar u_i(t)dy-d_i\bar u_i\leq0,
$$
and $u_i(0,x)\leq \bar u_i(0)=A$. A simple comparison principle asserts $(u_1,u_2)(t,x)\leq(\bar u_1,\bar u_2)(t)$ for $x\in[g(t),h(t)]$ and $t\geq0$.

Motivated by \cite[Theorem 3.3]{XLZ}, let us construct the iteration sequence
$\{(\bar u_1^{(n)},\bar u_2^{(n)})\}$ with initial value $(\bar u_1^{(0)},\bar u_2^{(0)}):=(A,A)$. It then follows from \cite{PP} that the limit of iteration sequences $\{(\bar u_1^{(n)},\bar u_2^{(n)})\}$ exists as $n\rightarrow\infty$, and we denote it by $(U_1^\Delta,U_2^\Delta)$, which satisfies
\begin{eqnarray*}
\left\{
\begin{array}{lll}
(U_1^\bigtriangleup)_t=b(t)U_2^\bigtriangleup-(a(t)+m_1(t)) U_1^\bigtriangleup-\alpha_1(t)(U_1^\bigtriangleup)2,  & t\in(0^+,\tau],\\[1mm]
(U_2^\bigtriangleup)_t=a(t)U_1^\bigtriangleup-m_2(t) U_2^\bigtriangleup-\alpha_2(t) (U_2^\bigtriangleup)^2, & t\in(0^+,\tau],\\[1mm]
U_1^\bigtriangleup(0^+)=U_1^\bigtriangleup(0),U_2^\bigtriangleup(0^+)=H(U_2^\bigtriangleup(0)),\\[1mm]
U_1^\bigtriangleup(0)=U_1^\bigtriangleup(\tau),U_2^\bigtriangleup(0)=U_2^\bigtriangleup(\tau) \\[1mm]
\end{array} \right.
\end{eqnarray*}
We further obtain
$$\lim\limits_{n\to+\infty}{(\bar u_1,\bar u_2)(t+n\tau)(t)}\leq(U_1^\Delta,U_2^\Delta)(t)$$
by the iteration method \cite[Theorem 4.5]{XLS}. Recalling $\underline{\lambda}(\infty,H'(0))\geq0$, $J_1(x)\neq J_2(x)$ and all coefficients are constant, one easily checks the $(U_1^\Delta,U_2^\Delta)=(0,0)$ similarly as Theorem \ref{th4.1}(i). Therefore, we get $\lim\limits_{t\to+\infty}{(u_1,u_2)(t,x)}=(0,0)$ uniformly in $x\in[g(t),h(t)]$.

\begin{cor}\label{co4.1}
Suppose $J_1(x)=J_2(x)$ and the principal eigenvalue satisfies $\lambda^*(\infty,H'(0))\geq0$, then
$$\lim\limits_{t\to+\infty}{(u_1,u_2)(t,x)}=(0,0)$$
uniformly in $x\in[g(t),h(t)]$.
\end{cor}

\vspace{3mm}
Before discussing the following part, we first give a conclusion for convenience.
\begin{prop}\label{prop4.5}
If $\lambda^*(b(t),a(t)+m_1(t),a(t),m_2(t);\infty;H'(0))<0$, then for any given $A^* \geq B:=\{\frac{b^M}{a^m+m_1^m},\frac{a^M}{m_2^m}\}$, the following problem
\begin{eqnarray*}
\left\{
\begin{array}{lll}
u_{1t}=b(t)u_2-(a(t)+m_1(t))u_1,  & t\in((n\tau)^+,(n+1)\tau],\\[1mm]
u_{2t}=a(t)u_1-m_2(t)u_2, & t\in((n\tau)^+,(n+1)\tau],\\[1mm]
u_1((n\tau)^+)=u_1(n\tau),\,u_2((n\tau)^+)=H(u_2(n\tau)),  \\[1mm]
u_i(0)=A^*, &i=1,2
\end{array} \right.
\end{eqnarray*}
admits a unique positive solution $(\bar u_1,\bar u_2)(t)$ and $\lim\limits_{n\to+\infty}{(\bar u_1,\bar u_2)(t+n\tau)(t)}=(u_1^*,u_2^*)(t)$ uniformly in $[0,\tau]$, where $(u_1^*,u_2^*)(t)$ satisfies
\begin{eqnarray}
\left\{
\begin{aligned}
&(u_1^*)_{t}=b(t)u_2^*-(a(t)+m_1(t))u_1^*,  && t\in(0^+,\tau],\\[1mm]
&(u_2^*)_t=a(t)u_1^*-m_2(t)u_2^*, && t\in(0^+,\tau],\\[1mm]
&u_1^*(0^+)=u_1^*(0),\,u_2^*(0^+)=H(u_2^*(0)),  \\[1mm]
&u_1^*(0)=u_1^*(\tau),\,u_2^*(0)=u_2^*(\tau).\\[1mm]
\end{aligned} \right.
\label{t1}
\end{eqnarray}
\end{prop}

\begin{lem}\label{le4.6}
Let $(u_1,u_2;(g,h))$ be the unique solution to \eqref{a01} with $J_1(x)=J_2(x)$. Assume that $\lambda^*(b(t),a(t)+m_1(t)+\alpha_1(t)u_1^*,a(t),m_2(t)+\alpha_2(t)u_2^*;\infty;H'(0))<0$ with $(u_1^*,u_2^*)$ defined in \eqref{t1}, if $h_\infty-g_\infty=\infty$, then $h_\infty=-g_\infty=\infty$.
\end{lem}
\bpf
Since $\lambda^*(b(t),a(t)+m_1(t)+\alpha_1(t)u_1^*,a(t),m_2(t)+\alpha_2(t)u_2^*;\infty;H'(0))<0$, we get $\lambda^*(b(t),a(t)+m_1(t),a(t),m_2(t);\infty;H'(0))<0$. It then follows from Proposition \ref{prop4.5} that the following initial value problem
\begin{eqnarray*}
\left\{
\begin{array}{lll}
u_{1t}=b(t)u_2-(a(t)+m_1(t))u_1,  & t\in((n\tau)^+,(n+1)\tau],\\[1mm]
u_{2t}=a(t)u_1-m_2(t)u_2, & t\in((n\tau)^+,(n+1)\tau],\\[1mm]
u_1((n\tau)^+)=u_1(n\tau),\,u_2((n\tau)^+)=H(u_2(n\tau),  \\[1mm]
u_i(0)=\overline A, &i=1,2\\[1mm]
\end{array} \right.
\end{eqnarray*}
admits a unique positive solution $(\bar u_1,\bar u_2)(t)$ and $\lim\limits_{n\to+\infty}{(\bar u_1,\bar u_2)(t+n\tau)(t)}=(u_1^*,u_2^*)(t)$ for $t\in[0,\tau]$, where $\overline A\geq\max\{A,A^*\}$ with $A$ defined in Lemma \ref{m2} and $A^*$ defined in Proposition \ref{prop4.5}. Also, a comparison principle yields $(u_1,u_2)(t,x)\leq (\bar u_1,\bar u_2)(t)$ for $t\geq0$ and $x\in(-\infty,+\infty)$.
So for any given $\epsilon_1$, there exists a large integer $N_1$ such that $(u_1(t,x),u_2(t,x))\leq(u_1^*(t)+\epsilon_1,u_2^*(t)+\epsilon_1)$ for $t\geq N_1\tau$ and $x\in[g(t),h(t)]$, which together with \eqref{a01} yields
\begin{eqnarray*}
\left\{
\begin{array}{lll}
u_{1t}-d_1\mathcal{L}[u_1]\geq b(t)u_2\\[1mm]
\qquad-[a(t)+m_1(t)+\alpha_1(t)(u_1^*+\epsilon_1)]u_1,  & t\in((n\tau)^+,(n+1)\tau],\,x \in (g(t),h(t)),\\[1mm]
u_{2t}-d_2\mathcal{L}[u_2]\geq a(t)u_1\\[1mm]
\qquad-[m_2(t)+\alpha_2(t)(u_2^*+\epsilon_1)]u_2, & t\in((n\tau)^+,(n+1)\tau],\,x \in (g(t),h(t)),\\[1mm]
u_1((n\tau)^+,x)=u_1(n\tau,x), & x \in (g(n\tau),h(n\tau)), \\[1mm]
u_2((n\tau)^+,x)=H(u_2(n\tau,x)), & x \in (g(n\tau),h(n\tau)), \\[1mm]
h'(t)=\sum\limits_{i=1}^{2}\mu_i\int_{g(t)}^{h(t)}\int_{h(t)}^{+\infty}J_i(x-y)u_i(t,x)dydx, &  t\in(n\tau,(n+1)\tau], \\[1mm]
g'(t)=-\sum\limits_{i=1}^{2}\mu_i\int_{g(t)}^{h(t)}\int_{-\infty}^{g(t)}J_i(x-y)u_i(t,x)dydx, &  t\in(n\tau,(n+1)\tau], \\[1mm]
u_1(t,x)=u_2(t,x)=0, & t\geq N_1\tau,\,x \in \{g(t),h(t)\}, \\[1mm]
n=N_1,N_1+1,N_1+2,\dots.\\[1mm]
\end{array} \right.
\end{eqnarray*}
Recalling $\lambda^*(b(t),a(t)+m_1(t)+\alpha_1(t)u_1^*,a(t),m_2(t)+\alpha_2(t)u_2^*;\infty;H'(0))<0$, so $\lambda^*(b(t),a(t)+m_1(t)+\alpha_1(t)(u_1^*+\epsilon_1),a(t),m_2(t)+\alpha_2(t)(u_2^*+\epsilon_1);\infty;H'(0))<0$ provided that positive constant $\epsilon_1$ is suitable small. Let $(\underline u_1,\underline u_2;(\underline g,\underline h))$ be the unique solution to the following problem
\begin{eqnarray*}
\left\{
\begin{array}{lll}
\underline u_{1t}-d_1\underline{\mathcal{L}}[\underline u_1]= b(t)\underline u_2\\[1mm]
\qquad-[a(t)+m_1(t)+\alpha_1(t)(u_1^*+\epsilon_1)]\underline u_1,  & t\in((n\tau)^+,(n+1)\tau],\,x \in (\underline g(t),\underline h(t)),\\[1mm]
\underline u_{2t}-d_2\underline{\mathcal{L}}[\underline u_2]= a(t)\underline u_1\\[1mm]
\qquad-[m_2(t)+\alpha_2(t)(u_2^*+\epsilon_1)]\underline u_2, & t\in((n\tau)^+,(n+1)\tau],\,x \in (\underline g(t),\underline h(t)),\\[1mm]
\underline u_1((n\tau)^+,x)=\underline u_1(n\tau,x), & x \in (\underline g(n\tau),\underline h(n\tau)), \\[1mm]
\underline u_2((n\tau)^+,x)=H(\underline u_2(n\tau,x)), & x \in (g(n\tau),h(n\tau)), \\[1mm]
\underline h'(t)=\sum\limits_{i=1}^{2}\mu_i\int_{\underline g(t)}^{\underline h(t)}\int_{\underline h(t)}^{+\infty}J_i(x-y)\underline u_i(t,x)dydx, &  t\in(n\tau,(n+1)\tau], \\[1mm]
\underline g'(t)=-\sum\limits_{i=1}^{2}\mu_i\int_{\underline g(t)}^{\underline h(t)}\int_{-\infty}^{\underline g(t)}J_i(x-y)\underline u_i(t,x)dydx, &  t\in(n\tau,(n+1)\tau], \\[1mm]
\underline u_1(t,x)=\underline u_2(t,x)=0, & t\geq N_1\tau,\,x \in \{\underline g(t),\underline h(t)\}, \\[1mm]
\underline u_i(N_1\tau,x)=u_i(N_1\tau,x), & x\in [\underline g(N_1\tau),\underline h(N_1\tau)],\\[1mm]
\underline h(N_1\tau)=-\underline g(N_1\tau)=h(N_1\tau),\\[1mm]
i=1,2,\,n=N_1,N_1+1,N_1+2,\dots,\\[1mm]
\end{array} \right.
\end{eqnarray*}
where $\underline{\mathcal{L}}[\underline u_i]=\int_{\underline g(t)}^{\underline h(t)}J(x-y)\underline u_i(t,y)dy-\underline u_i(t,x)$. It follows from the comparison principle that
$$
(u_1,u_2)(t,x)\geq(\underline u_1,\underline u_2)(t,x),\,\,h(t)\geq\underline h(t),\,\,\underline g(t)\geq g(t),\,\,t\geq N_1\tau,\,x\in(\underline g(t),\underline h(t)).
$$

In the following, we will prove $\underline h_\infty-\underline g_\infty=\infty$. On the contrary, we have $\underline h_\infty-\underline g_\infty<\infty$, it then follows from Lemma \ref{le4.3} that $$\lambda^*(b(t),a(t)+m_1(t)+\alpha_1(t)(u_1^*+\epsilon_1),a(t),m_2(t)+\alpha_2(t)(u_2^*+\epsilon_1);(\underline g_\infty,\underline h_\infty);H'(0))\geq0.$$ By
the monotonicity in Lemma 3.2 (iii), one easily checks that $$\lambda^*(b(t),a(t)+m_1(t)+\alpha_1(t)(u_1^*+\epsilon_1),a(t),m_2(t)+\alpha_2(t)(u_2^*+\epsilon_1);(\underline g(N_1\tau),\underline h(N_1\tau));H'(0))>0,$$
which together with $\lambda^*(b(t),a(t)+m_1(t)+\alpha_1(t)(u_1^*+\epsilon_1),a(t),m_2(t)+\alpha_2(t)(u_2^*+\epsilon_1);\infty;H'(0))<0$, yields
$$\lambda^*(b(t),a(t)+m_1(t)+\alpha_1(t)(u_1^*+\epsilon_1),a(t),m_2(t)+\alpha_2(t)(u_2^*+\epsilon_1);(-L^*,L^*);H'(0))=0$$ for some positive constant $L^*$. Since $h_\infty-g_\infty=\infty$, we can enlarge the integer $N_1$ such that $h(N_1\tau)-g(N_1\tau)>2L^*$, so $\underline h(N_1\tau)-\underline g(N_1\tau)>2L^*$. Therefore, by Lemma \ref{le3.2} (i) and (iii), we obtain
$$\begin{array}{llllll}
&&\lambda^*(b(t),a(t)+m_1(t)+\alpha_1(t)(u_1^*+\epsilon_1),a(t),m_2(t)+\alpha_2(t)(u_2^*+\epsilon_1);(\underline g_\infty,\underline h_\infty);H'(0))\\[1mm]
&\leq&\lambda^*(b(t),a(t)+m_1(t)+\alpha_1(t)(u_1^*+\epsilon_1),a(t),m_2(t)+\alpha_2(t)(u_2^*+\epsilon_1);(\underline g(N_1\tau),\underline h(N_1\tau));H'(0))\\[1mm]
&<&\lambda^*(b(t),a(t)+m_1(t)+\alpha_1(t)(u_1^*+\epsilon_1),a(t),m_2(t)+\alpha_2(t)(u_2^*+\epsilon_1);(-L^*,L^*);H'(0))\\[1mm]
&=&0,
\end{array}$$
which leads a contradiction. So $\underline h_\infty-\underline g_\infty=\infty$.

The proof ends with $\underline h_\infty=-\underline g_\infty=\infty$. Once it has done, this together with $h(t)\geq\underline h(t)$ and $\underline g(t)\geq g(t)$, to obtain $h_\infty=-g_\infty=\infty$. Without loss of generality, suppose on the contrary that $\underline h_\infty<\infty$ and $-\underline g_\infty=\infty$. Since $\lambda^*(b(t),a(t)+m_1(t)+\alpha_1(t)(u_1^*+\epsilon_1),a(t),m_2(t)+\alpha_2(t)(u_2^*+\epsilon_1);\infty;H'(0))<0$, for any sufficiently small $\epsilon_2$, there exists sufficiently large $g_1$ such that $\lambda^*(b(t),a(t)+m_1(t)+\alpha_1(t)(u_1^*+\epsilon_1),a(t),m_2(t)+\alpha_2(t)(u_2^*+\epsilon_1);(-g_1,\underline h_\infty-\epsilon_2);H'(0))<0$. For such $g_1$ and $\epsilon_2$, we can select a large integer $N_2$ such that
$$
(-g_1,\underline h_\infty-\epsilon_2)\subset(\underline g(t),\underline h(t)),\forall t\geq N_2\tau.
$$
Let $(\widehat{u}_1,\widehat{u}_2;(-g_1,\underline h_\infty-\epsilon_2))$ be the unique solution to
\begin{eqnarray*}
\left\{
\begin{array}{lll}
\widehat u_{1t}-d_1\widehat{\mathcal{L}}[\widehat u_1]= b(t)\widehat u_2\\[1mm]
\qquad-[a(t)+m_1(t)+\alpha_1(t)(u_1^*+\epsilon_1)]\widehat u_1,  & t\in((n\tau)^+,(n+1)\tau],\,x \in (-g_1,\underline h_\infty-\epsilon_2),\\[1mm]
\widehat u_{2t}-d_2\widehat{\mathcal{L}}[\widehat u_2]= a(t)\widehat u_1\\[1mm]
\qquad-[m_2(t)+\alpha_2(t)(u_2^*+\epsilon_1)]\widehat u_2, & t\in((n\tau)^+,(n+1)\tau],\,x \in (-g_1,\underline h_\infty-\epsilon_2),\\[1mm]
\widehat u_1((n\tau)^+,x)=\widehat u_1(n\tau,x), & x \in (-g_1,\underline h_\infty+\epsilon_2), \\[1mm]
\widehat u_2((n\tau)^+,x)=H(\widehat u_2(n\tau,x)), & x \in (-g_1,\underline h_\infty+\epsilon_2), \\[1mm]
\widehat u_i(N_2\tau,x)=\underline u_i(N_2\tau,x),& x\in [-g_1,\underline h_\infty+\epsilon_2],\\[1mm]
i=1,2,\,n=N_2,N_2+1,N_2+2,\dots,\\[1mm]
\end{array} \right.
\end{eqnarray*}
where $\widehat{\mathcal{L}}[u]=\int_{-g_1}^{\underline h_\infty+\epsilon_2}J(x-y)u(t,y)dy-u(t,x)$. The comparison principle yields
$$
(\underline u_1,\underline u_2)(t,x)\geq(\widehat{u}_1,\widehat{u}_2)(t,x),\, \forall t\geq N_2\tau, \,x\in[-g_1,\underline h_\infty+\epsilon_2].
$$
Analogously as methods in Lemma \ref{le4.3}, we can find a constant $\delta^*$ and a large integer $N_3(>N_2)$ such that
$\underline h'(t+n\tau)\geq \delta^*>0$ for $t\geq N_3\tau$. So $\underline h(t+n\tau)\rightarrow \infty$ as $t\rightarrow \infty$, which leads a contradiction to $\underline h_\infty<\infty$. The proof is now complete.
\epf

\begin{lem}\label{le4.7}
Assume $\underline{\lambda}((-h_0,h_0),H'(0))>0$ and $\|u_{1,0}(x)\|_{C([-h_0,h_0])}+\|u_{2,0}(x)\|_{C([-h_0,h_0])}$ is sufficiently small, then $h_\infty-g_\infty<\infty$ and $\lim\limits_{t\to+\infty}{||u_1(t,\cdot)||_{C[g(t),h(t)]}}=\lim\limits_{t\to+\infty}{||u_2(t,\cdot)||_{C[g(t),h(t)]}}=(0,0)$ uniformly for $x\in [g(t),h(t)]$.
\end{lem}
\bpf
We see from $\underline{\lambda}((-h_0,h_0),H'(0))>0$ that $\underline{\lambda}((-h_1,h_1),H'(0))>0$ for some $h_1=h_0+\epsilon$ with $\epsilon>0$ small. Let $(\phi(t,x),\psi(t,x))$ be the corresponding normalized eigenfunction pair in the fixed domain $[-h_1,h_1]$ and $||\phi(t,x)||_{C[-h_1,h_1]},\,||\psi(t,x)||_{C[-h_1,h_1]}\leq1$ for $t\in[0,\tau]$.

Define
$$
\overline W_1(t,x)=C_1e^{-\gamma t}\phi_{h_1}(t,x),\,\,\overline W_2(t,x)=C_1e^{-\gamma t}\psi_{h_1}(t,x),\,\,t\geq0,\,\,x\in[-h_1,h_1],
$$
$$
\overline \eta(t)=h_0+(h_1-h_0)(1-e^{-\gamma t}),\,\,t\geq0,
$$
where positive constants $\gamma$ and $C_1$ to be chosen later. It is obvious that $\overline \eta(t)\in[h_0,h_1)$.

Careful calculations in $t\in((n\tau)^+,(n+1)\tau]$ and $x\in (-\overline \eta(t),\overline \eta(t))$ yield
$$\begin{array}{llllll}
&&\overline W_{1t}-d_1[\int_{-\overline \eta(t)}^{\overline \eta(t)}J_1(x-y)\overline{W}_1(t,y)dy-\overline{W}_1(t,x)]-b(t)\overline W_2+[a(t)+m_1(t)]\overline W_1+\alpha_1(t)\overline W_1^2\\[1mm]
&=&C_1e^{-\gamma t}[-\gamma\phi_{h_1}+(\phi_{h_1})_t-d_1(\int_{-\overline \eta(t)}^{\overline \eta(t)}J_1(x-y)\phi_{h_1}(t,y)dy-\phi_{h_1})\\[1mm]
&-&b(t)\psi_{h_1}+[a(t)+m_1(t)]\phi_{h_1}+C_1\alpha_1(t)e^{-\gamma t}(\phi_{h_1})^2]\\[1mm]
&\geq&C_1e^{-\gamma t}[-\gamma\phi_{h_1}+(\phi_{h_1})_t-d_1(\int_{-h_1}^{h_1}J_1(x-y)\phi_{h_1}(t,y)dy-\phi_{h_1})\\[1mm]
&-&b(t)\psi_{h_1}+(a(t)+m_1(t))\phi_{h_1}]\\[1mm]
&\geq&C_1e^{-\gamma t}\phi_{h_1}[\underline{\lambda}((-h_1,h_1),H'(0))-\gamma]\\[1mm]
&>0&
\end{array}$$
provided that $\gamma:=\frac{\underline{\lambda}((-h_1,h_1),H'(0))}{2}>0$.

Similarly, it is clear that
$$\begin{array}{llllll}
&&\overline W_{2t}-d_2(\int_{-\overline \eta(t)}^{\overline \eta(t)}J_1(x-y)\overline{W}_2(t,y)dy-\overline{W}_2(t,x))-a(t)\overline W_1+m_2(t)\overline W_2+\alpha_2(t)\overline W_2^2\\[1mm]
&=&C_1e^{-\gamma t}[-\gamma\psi_{h_1}+(\psi_{h_1})_t-d_2(\int_{-\overline \eta(t)}^{\overline \eta(t)}J_2(x-y)\psi_{h_1}(t,y)dy-\psi_{h_1})\\[1mm]
&-&a(t)\phi_{h_1}+m_2(t)\psi_{h_1}+C_1\alpha_2(t)e^{-\gamma t}(\psi_{h_1})^2]\\[1mm]
&\geq&C_1e^{-\gamma t}\psi_{h_1}[\underline{\lambda}((-h,h),H'(0))-\gamma]\\[1mm]
&>&0.
\end{array}$$
Meanwhile, for $x\in[-\overline\eta(n\tau),\overline\eta(n\tau)]$, we get $$\overline W_1((n\tau)^+,x)=C_1e^{-\gamma n\tau}\phi_{h_1}((n\tau)^+,x)=\overline W_1(n\tau,x)$$ and
$$\begin{array}{llllll}
&&\overline W_2((n\tau)^+,x)=C_1e^{-\gamma n\tau}\psi_{h_1}((n\tau)^+,x)\\[1mm]
&=&C_1e^{-\gamma n\tau}H'(0)\psi_{h_1}(n\tau,x)\\[1mm]
&=&H'(0)\overline W_2(n\tau,x)\\[1mm]
&\geq&H(\overline W_2(n\tau,x))
\end{array}$$
according to ($\mathcal{A}$).

Since $[-\overline \eta(t),\overline \eta(t)]\subset(-h_1,h_1)$, we deduce that
$$\begin{array}{llllll}
&&\sum\limits_{i=1}^{2}\mu_i\int_{-\overline \eta(t)}^{\overline \eta(t)}\int_{\overline \eta(t)}^{+\infty}J_i(x-y)\overline W_i(t,x)dydx\\[1mm]
&\leq&\sum\limits_{i=1}^{2}\mu_i\int_{-\overline \eta(t)}^{\overline \eta(t)}\overline W_i(t,x)dx\\[1mm]
&<&C_1e^{-\gamma t}\int_{-h_1}^{h_1}(\mu_1\phi_{h_1}(t,x)+\int_{-h_1}^{h_1}\mu_2\psi_{h_1}(t,x))dy\\[1mm]
&\leq&2h_1C_1e^{-\gamma t}(\mu_1+\mu_2)\\[1mm]
\end{array}$$
and
$$
\overline\eta'(t)=(h_1-h_0)\gamma e^{-\gamma t},
$$
so $$\overline\eta'(t)\geq\sum\limits_{i=1}^{2}\mu_i\int_{-\overline \eta(t)}^{\overline \eta(t)}\int_{\overline \eta(t)}^{+\infty}J_i(x-y)\overline W_i(t,x)dydx.$$
provided that $C_1:=(h_1-h_0)\frac{\gamma}{2h_1(\mu_1+\mu_2)}>0$.

Recalling $[-\overline \eta(t),\overline \eta(t)]\subset(-h_1,h_1)$, we obtain $\overline W_1(t,\pm\overline \eta(t))=C_1e^{-\gamma t}\phi_{h_1}(t,\pm\overline \eta(t))>0$ and $\overline W_2(t,\pm\overline \eta(t))=C_1e^{-\gamma t}\psi_{h_1}(t,\pm\overline \eta(t))>0$ hold for $i=1,2$ and $t>0$.
Also,
$$\overline W_1(0,x)=C_1\phi_{h_1}(0,x)\geq u_{1,0}(x),\,\overline W_2(0,x)=C_1\psi_{h_1}(0,x)\geq u_{2,0}(x),\, x\in[-h_0,h_0]$$
provided that
$$\|u_{1,0}(x)\|_{C([-h_0,h_0])}+\|u_{2,0}(x)\|_{C([-h_0,h_0])}\leq C_1\min\{\min_{x\in[-h_0,h_0]}\psi_{h_1}(0,x),\min_{x\in[-h_0,h_0]}\psi_{h_1}(0,x)\}.$$

Therefore, $(\overline W_1,\overline W_2;(-\overline \eta,\overline\eta))$ is an upper solution to problem \eqref{a01}. It follows from Lemma \ref{le4.2} that $(u_1(t,x),u_2(t,x))\leq (\overline W_1,\overline W_2)(t,x)$ in $(t,x)\in(0,+\infty)\times(g(t),h(t))$, and $h(t)\leq \overline\eta(t)$, $g(t)\geq -\overline\eta(t)$ for $t\in(0,+\infty)$. Since $\lim\limits_{t\to+\infty}{(\overline W_1,\overline W_2)(t,x)}=(0,0)$ in $x\in[g(t),h(t)]$, we derive that $\lim\limits_{t\to+\infty}{||(u_1(t,\cdot),u_2(t,\cdot))||_{C[g(t),h(t)]}}=(0,0).$ Meanwhile, $h_\infty-g_\infty=\lim\limits_{t\to+\infty}{(h(t)-g(t))}\leq 2h_1<+\infty$ can be obtained. This finishes the proof.
\epf

\begin{lem}\label{le4.8}
If $-g_\infty=h_\infty=\infty$ and $\overline{\lambda}(\infty,H'(0))<0$, then $$\lim\limits_{m\to+\infty}{(u_1,u_2)(t+m\tau,x)}=(U^\triangle,V^\triangle)(t)$$
uniformly for $t\in[0,\tau]$ and locally uniformly for $x\in(-\infty,+\infty)$, where $(U^\triangle,V^\triangle)(t)$ is the unique positive solution to problem
\begin{eqnarray}
\left\{
\begin{array}{lll}
U_t=b(t)V-(a(t)+m_1(t))U-\alpha_1(t)U^2,  & t\in(0^+,\tau],\\[1mm]
V_t=a(t)U-m_2(t) V-\alpha_2(t) V^2, & t\in(0^+,\tau],\\[1mm]
U(0^+)=U(0),V(0^+)=H(V(0)),\\[1mm]
U(0)=U(\tau),V(0)=V(\tau). \\[1mm]
\end{array} \right.
\label{c20}
\end{eqnarray}
\end{lem}
\bpf It is clear that the solution $(U^\triangle,V^\triangle)(t)$ of spatial-independent problem \eqref{c20} exists uniquely \cite[Theorem 3.6]{XLZ} and we omit the details here.

In the following, we claim that
\begin{equation}
\begin{array}{lll}
\liminf_{m\to+\infty}{(u_1,u_2)(t+m\tau,x)}\geq(U^\triangle,V^\triangle)(t)\text{\,\,locally uniformly in\,\,}[0,\tau]\times(-\infty,+\infty)
\end{array}
\label{c200}
\end{equation}
and
\begin{equation}
\begin{array}{lll}
\limsup_{m\to+\infty}{(u_1,u_2)(t+m\tau,x)}\leq(U^\triangle,V^\triangle)(t)\text{\,\,uniformly in\,\,}[0,\tau]\times(-\infty,+\infty).
\end{array}
\label{c201}
\end{equation}

We first prove that \eqref{c200} holds. Since $\overline{\lambda}(\infty,H'(0))<0$, there exists a positive constant $L$ such that $$\overline{\lambda}((-L,L),H'(0))<0.$$ Also, by $-g_\infty=h_\infty=\infty$, we can find the positive integer $n_L$ such that $h(t)-g(t)\geq 2L$ for any $t\geq n_L\tau$.

We first consider the following initial boundary problem
\begin{eqnarray}
\left\{
\begin{array}{lll}
\underline{u}_{1t}-d_1\mathcal{\tilde{L}}_1[\underline u_1]=b(t)\underline u_2-(a(t)+m_1(t))\underline u_1\\
\qquad\qquad\qquad\quad-\alpha_1(t)\underline u_1^2,  & t\in((n\tau)^+,(n+1)\tau],\,x \in [-L,L],\\[1mm]
\underline u_{2t}-d_2\mathcal{\tilde{L}}_2[\underline u_2]=a(t)\underline u_1-m_2(t)\underline u_2-\alpha_2(t)\underline u_2^2, & t\in((n\tau)^+,(n+1)\tau],\,x \in [-L,L],\\[1mm]
\underline u_1((n\tau)^+,x)=\underline u_1(n\tau,x),&x \in (-L,L),\,n\geq n_L,\\[1mm]
\underline u_2((n\tau)^+,x)=H(\underline u_2(n\tau,x)), &x \in (-L,L),\,n\geq n_L,\\[1mm]
\underline u_1(n_L\tau,x)=u_1(n_L\tau,x),\,\underline u_2(n_L\tau,x)=u_2(n_L\tau,x),& x\in[-L,L],\,n=n_L,n_L+1,\dots, \\[1mm]
\end{array} \right.
\label{c21}
\end{eqnarray}
in a fixed domain $[-L,L]$. A comparison principle to conclude
$$(\underline u_1,\underline u_2)(t,x)\leq(u_1,u_2)(t,x),\,\,\forall t\geq n_L\tau,\,\,-L\leq x\leq L.$$
Since $\overline{\lambda}(\infty,H'(0))<0$, by Theorem \ref{th4.1} (iii), problem \eqref{c21} admits a unique positive steady state $(U(t,x),V(t,x))$ defined   in problem \eqref{d02} with $[L_1,L_2]$ replaced by $[-L,L]$, and
$$
\liminf_{m\to+\infty}{(u_1,u_2)(t+m\tau,x)}\geq(U,V)(t,x)\,\,\text{uniformly\,\,on\,\,} [0,\tau]\times[-L,L],
$$
Next letting $L\rightarrow\infty$ yileds
$$
\liminf_{m\to+\infty}{(u_1,u_2)(t+m\tau,x)}\geq (U^\triangle,V^\triangle)(t)\,\,\text{locally\,\,uniformly\,\,in\,\, }[0,\tau]\times(-\infty,+\infty).
$$

Finally we prove \eqref{c201} holds. Let $(\bar u_1,\bar u_2)$ be a solution to the following problem
\begin{eqnarray}
\left\{
\begin{array}{lll}
u_{1t}=b(t)u_2-(a(t)+m_1(t))u_1-\alpha_1(t) u_1^2,  & t\in((n\tau)^+,(n+1)\tau],\\[1mm]
u_{2t}=a(t) u_1-m_2(t) u_2-\alpha_2(t) u_2^2, & t\in((n\tau)^+,(n+1)\tau],\\[1mm]
u_1((n\tau)^+)= u_1(n\tau),\, u_2((n\tau)^+)=H( u_2(n\tau)),\\[1mm]
u_1(0)=u_2(0)=A,\\[1mm]
\end{array} \right.
\label{c22}
\end{eqnarray}
where $A$ is defined in Lemma \ref{m2}. Recalling the diffusion term in \eqref{a01} and
$$
-d_i(\int_{g(t)}^{h(t)}J_i(x-y)\bar u_i(t)dy-\bar u_i(t))\geq 0,\,\,i=1,2,
$$
which together with a comparison principle yields $(u_1,u_2)(t,x)\leq(\bar u_1,\bar u_2)(t)$ for $t\in[0,+\infty)$ and $x\in[g(t),h(t)]$, thus
\begin{equation}
\limsup_{m\to+\infty}{(u_1,u_2)(t+m\tau,x)}\leq\lim\limits_{m\to+\infty}{(\bar u_1,\bar u_2)(t+m\tau)}
\label{c23}
\end{equation}
uniformly for $t\in[0,+\infty)$ and $x\in[g(t+m\tau),h(t+m\tau)]$.

In the following, We construct iteration sequences $\{\tilde U^ {(m)}\}$ and $\{\tilde V^ {(m)}\}$ satisfying
\begin{eqnarray*}
\left\{
\begin{array}{lll}
\tilde U_t^{(m)}+K^*\tilde U^{(m)}=K^* \tilde U^{(m-1)}+b(t)\tilde V^{(m-1)}-(a(t)+m_1(t))\tilde U^{(m-1)}\\[1mm]
\qquad \qquad \quad \quad\quad-\alpha_1(t) (\tilde U^{(m-1)})^2,  & t\in(0^+,\tau],\\[1mm]
\tilde V_t^{(m)}+K^*\tilde V^{(m)}=K^* \tilde V^{(m-1)}+a(t)\tilde U^{(m-1)}-m_2(t)\tilde V^{(m-1)}\\[1mm]
\qquad \qquad \quad \quad\quad-\alpha_2(t)(\tilde V^{(m-1)})^2, & t\in(0^+,\tau],\\[1mm]
\tilde U^{(m)}(0^+)=\tilde U^{(m-1)}(\tau),\\[1mm]
\tilde V^{(m)}(0^+)=H(\tilde V^{(m-1)}(\tau)),\\[1mm]
\tilde U^{(m)}(0)=\tilde U^{(m-1)}(\tau),\,\tilde V^{(m)}(0)=\tilde V^{(m-1)}(\tau).\\[1mm]
\end{array} \right.
\end{eqnarray*}
Similarly  as Theorem \ref{th4.1}, the induction for $m$ yields
\begin{equation}
(\bar u_1,\bar u_2)(t+m\tau)\leq(\tilde U^ {(m)},\tilde V^ {(m)})(t),\,t\in[0,\tau],\,m=0,1,2,\dots.
\label{c24}
\end{equation}
Recalling that problem \eqref{c20} has a unique solution $(U^\triangle, V^\triangle)$, then iteration sequence satisfies
\begin{equation}
\lim\limits_{m\to+\infty}{(\tilde U^ {(m)},\tilde V^ {(m)})(t)}=(U^\triangle, V^\triangle)(t),\,t\in[0,\tau].
\label{c25}
\end{equation}
It follows from \eqref{c23}, \eqref{c24} and \eqref{c25} that $\limsup_{m\to+\infty}{(u_1,u_2)(t+m\tau,x)}\leq(U^\triangle,V^\triangle)(t)$
uniformly in $[0,\tau]\times(-\infty,+\infty)$.
\epf


\vspace{3mm}

Usually, we say vanishing if $\lim\limits_{t\to+\infty}{(u_1,u_2)(t,x)}=(0,0)$ in $x\in[g(t),h(t)]$, and spreading if $h_\infty-g_\infty=\infty$ and $\liminf_{m\to+\infty}{(u_1,u_2)(t+m\tau,x)}>0$ in $[0,\tau]\times[g(t),h(t)]$. Recalling Lemmas \ref{le4.3}, \ref{le4.4} and \ref{le4.8}, we have the following sufficient conditions for spreading and vanishing.

\begin{thm}\label{th4.9} Assume that $\mathbf{(J)}$ holds. We have the following statements:\\
(i) If $\underline{\lambda}(\infty,H'(0))\geq0$, $J_1(x)\neq J_2(x)$ and all coefficients are constant, then vanishing occurs;\\
(ii) If $\overline{\lambda}(\infty,H'(0))<0$ and

(a) $-g_\infty=h_\infty=\infty$, then spreading happens;

(b) $h_\infty-g_\infty<\infty$, $J_1(x)=J_2(x)$, then vanishing occurs.
\end{thm}

Specially, if the principle eigenvalue exists, we have the following spreading-vanishing dichotomy.
\begin{thm}\label{th4.10} Suppose $J_1(x)=J_2(x)$ holds, that is, the principal eigenvalue $\lambda^*$ exists, then the following assertions hold:\\
(i) If $\lambda^*(\infty,H'(0))\geq0$, then vanishing occurs;\\
(ii) If $\lambda^*(\infty,H'(0))<0$ and

(a) $-g_\infty=h_\infty=\infty$, then spreading happens;

(b) $h_\infty-g_\infty<\infty$, then vanishing occurs.
\end{thm}

In the following, the spreading-vanishing dichotomy in expanding capacities $(\mu_i)(i=1,2)$ are investigated.

\begin{thm}\label{th4.11} Assume that $\mathbf{(J)}$ holds and $J_1(x)=J_2(x)$. \\
$(i)$ If $\lambda^*(\infty,H'(0))\geq0$, then vanishing occurs for any given $\mu_1>0$, $\mu_2>0$ and initial value $(u_{1,0},u_{2,0})$ satisfying \eqref{a02};\\
$(ii)$ If $\lambda^*((-h_0,h_0),H'(0))\leq0$, then $-g_\infty=h_\infty=\infty$ and spreading occurs for any given $\mu_1$, $\mu_2$ and initial value;\\
$(iii)$ If $\lambda^*((-h_0,h_0),H'(0))>0$ and $\lambda^*(\infty,H'(0))<0$, then there exists $\mu^{**}\geq\mu_{**}>0$ such that vanishing happens for $0<\mu_1+\mu_2\leq\mu_{**}$ and spreading happens for $\mu_1+\mu_2>\mu^{**}$ for any given initial datum;\\
$(iv)$ If initial value $(u_{1,0},u_{2,0})$ satisfying \eqref{a02} is sufficiently small, then vanishing occurs for any fixed $\mu_i(i=1,2)>0$.
\end{thm}
\bpf
$(i)$ If $\lambda^*(\infty,H'(0))\geq0$, it indicates by Corollary 4.1 that $(u_1,u_2)(t,x)\rightarrow(0,0)$ in $x\in[g(t),h(t)]$ as $t\rightarrow\infty$. Therefore vanishing occurs for any $\mu_1$, $\mu_2$ and initial value satisfying \eqref{a02}.

$(ii)$ If $\lambda^*((-h_0,h_0),H'(0))\leq0$, we claim that $-g_\infty=h_\infty=\infty$ and spreading occurs. Otherwise, we first suppose $h_\infty-g_\infty<\infty$. It is derived from Lemma 4.3 that $\lambda^*((g_\infty,h_\infty),H'(0))\geq0$. So $\lambda^*((-h_0,h_0),H'(0))>\lambda^*((g_\infty,h_\infty),H'(0))\geq 0$ holds by
strictly monotonic decreasing of $\lambda^*$ with respect to $\Omega$ in Lemma \ref{le3.2} (iii), which contradicts to $\lambda^*((-h_0,h_0),H'(0))\leq0$, so $-g_\infty=h_\infty=\infty$. Also $\lambda^*((-h_0,h_0),H'(0))\leq0$ indicates $\lambda^*(\infty,H'(0))<0$, and by the same method in Lemma \ref{le4.8}, we prove that species spread.

$(iii)$ We first construct the upper solution, which is similar to Lemma \ref{le4.7} and we emphasize the difference and give the sketches here.
Modify
$$
\overline \eta(t)=h_0+2h_1(\mu_1+\mu_2)C_1\int_0^te^{-\gamma t}dt,\,\,t\geq0
$$
where $$\gamma:=\frac{\lambda^*((-h_1,h_1),H'(0))}{2}$$ is defined in Lemma \ref{le4.7} and
$$
C_1\geq\max\{\frac{\|u_{1,0}(x)\|_{C([-h_0,h_0])}+\|u_{2,0}(x)\|_{C([-h_0,h_0])}}{\min_{x\in[-h_0,h_0]}\phi_{h_1}(0,x)},
\frac{\|u_{1,0}(x)\|_{C([-h_0,h_0])}+\|u_{2,0}(x)\|_{C([-h_0,h_0])}}{\min_{x\in[-h_0,h_0]}\psi_{h_1}(0,x)}\}.
$$
It can be computed that
$$
\overline \eta(t)=h_0+2h_1(\mu_1+\mu_2)C_1\int_0^te^{-\gamma t}dt<h_0+2h_1(\mu_1+\mu_2)\frac{C_1}{\gamma}\leq h_1
$$
provided that
$$0<\mu_1+\mu_2\leq\frac{(h_1-h_0)\gamma}{2h_1C_1}=\mu_{**}.$$
The choice of $C_1$ satisfies $\overline W_1(0,x)=C_1\phi_{h_1}(0,x)\geq u_{1,0}(x)$ and $\overline W_2(0,x)=C_1\psi_{h_1}(0,x)\geq u_{2,0}(x)$ in $x\in[-h_0,h_0]$. Also,
$$\begin{array}{llllll}
&&\sum\limits_{i=1}^{2}\mu_i\int_{-\overline \eta(t)}^{\overline \eta(t)}\int_{\overline \eta(t)}^{+\infty}J_i(x-y)\overline W_i(t,x)dydx\\[1mm]
&\leq&2h_1(\mu_1+\mu_2)C_1e^{-\gamma t}\\[1mm]
&=&\bar \eta'(t).
\end{array}$$
Thus $(\overline W_1,\overline W_2;(-\overline \eta,\overline\eta))$ is an upper solution to problem \eqref{a01}, which indicates vanishing occurs if $0<\mu_1+\mu_2\leq\mu_{**}.$

Next, we claim that there exists a large $\mu^{**}$ such that spreading happens when $\mu_1+\mu_2>\mu^{**}$. In order to emphasize the dependence of the unique solution $(u_1,u_2;(g,h))$ of \eqref{a01} on $\mu$, we denote it by $(u_{1,\mu},u_{2,\mu};(g_\mu,h_\mu))$.

Recalling that $\lambda^*((-h_0,h_0),H'(0))>0$ and $\lambda^*(\infty,H'(0))<0$, it is easy to see that  $\lambda^*((g_{\mu^*}(t^*),h_{\mu^*}(t^*)),H'(0))=0$ for some $t^*,\mu^*>0$. Define $$2l^*:=h_{\mu^*}(t^*)-g_{\mu^*}(t^*),$$ it follows from the properties of $\lambda^*$ in Lemma \ref{le3.2} (iii) that $\lambda^*((g_{\mu_0}(t_0),h_{\mu_0}(t_0)),H'(0))<0$ provided that $h_{\mu_0}(t_0)-g_{\mu_0}(t_0)>2l^*$ for some given $t_0$ and $\mu_0$. So in the following, we firstly prove that there exists positive $\mu_0$, such that
\begin{equation}
h_{\mu_0}(t_0)-g_{\mu_0}(t_0)>2l^* \,\,\text{for\,\,some}\,\, t_0>0.
\label{c27}
\end{equation}
Suppose on the contrary that $h_{\mu}(t)-g_{\mu}(t)\leq2l^*$ for any $t,\mu>0$. Since that $-g_\mu(t)$ and $h_\mu(t)$ are strictly increasing in $t$ and $\mu$, so
$$
G_\infty:=\lim\limits_{t,\mu\rightarrow \infty}{g_\mu(t)}\,\,\text{and}\,\,H_\infty:=\lim\limits_{t,\mu\rightarrow \infty}{h_\mu(t)},
$$
moreover, $H_\infty-G_\infty\leq 2l^*$. Since that $J_i(0)>0$, there exists $\epsilon_0$ and $\kappa_0$ such that $J_i(x)>\kappa_0$ for $|x|<\epsilon_0$. Also, for such $\epsilon_0$, there exists positive constant $\mu_0$ and positive integer $N_0$ large enough such that
$$
h_\mu(t)+\epsilon_0/4>H_\infty \,\,\text{for}\,\,\mu\geq \mu_0\,\,\text{and}\,\,t\geq N_0\tau.
$$
Recalling equations in \eqref{a01} and integration in $t$ over $[N_0\tau,(N_0+1)\tau]$,  yields
$$\begin{array}{llllll}
&&h_\mu((N_0+1)\tau)-h_\mu(N_0\tau)=\sum\limits_{i=1}^{2}\mu_i\int_{N_0\tau}^{(N_0+1)\tau}
\int_{g_\mu(t)}^{h_\mu(t)}\int_{h_\mu(t)}^{+\infty}J_i(x-y)u_{i,\mu}(t,x)dydxdt\\[1mm]
&\geq&\sum\limits_{i=1}^{2}\mu_i\int_{N_0\tau}^{(N_0+1)\tau}
\int_{g_{\mu_0}(t)}^{h_{\mu_0}(t)}\int_{h_{\mu_0}(t)+\epsilon_0/4}^{+\infty}J_i(x-y)u_{i,\mu_0}(t,x)dydxdt\\[1mm]
&\geq&\sum\limits_{i=1}^{2}\mu_i\int_{N_0\tau}^{(N_0+1)\tau}
\int_{h_{\mu_0}(t)-\epsilon_0/2}^{h_{\mu_0}(t)}\int_{h_{\mu_0}(t)+\epsilon_0/4}^{h_{\mu_0}(t)+\epsilon_0/2}J_i(x-y)u_{i,\mu_0}(t,x)dydxdt\\[1mm]
&\geq&\epsilon_0\kappa_0/4\sum\limits_{i=1}^{2}\mu_i\int_{N_0\tau}^{(N_0+1)\tau}
\int_{h_{\mu_0}(t)-\epsilon_0/2}^{h_{\mu_0}(t)}(u_{i,\mu_0})_{\min}(t,x)dxdt,\\[1mm]
\end{array}$$
So $\mu_1+\mu_2\leq\frac{4[h_\mu((N_0+1)\tau)-h_\mu(N_0\tau)]}{\epsilon_0\kappa_0\sum\limits_{i=1}^{2}\mu_i\int_{N_0\tau}^{(N_0+1)\tau}
\int_{h_{\mu_0}(t)-\epsilon_0/2}^{h_{\mu_0}(t)}(u_{i,\mu_0})_{\min}(t,x)dxdt}<\infty$ (the right side of above equality has an upper bound that is independent of $\mu$ thanks to $h_{\mu}(t)-g_{\mu}(t)\leq2l^*$ for any $t,\mu>0$), which leads a contradiction. Therefore, \eqref{c27} holds and $\lambda^*((g_{\mu_0}(t_0),h_{\mu_0}(t_0)),H'(0))<0$ for some $t_0,\mu_0>0$. It follows from \cite[Lemma 4.11]{DNNNN} that spreading often happens for $\mu=\mu_0$. Since that $-g_{\mu}(t)$ and $h_{\mu}(t)$ are increasing in $\mu$, species spread for $\mu\geq\mu_0$, which ends the proof.

The conclusion $(iv)$ can be deduced by the method in Lemma \ref{le4.7} analogously.
\epf

\vspace{5mm}

We finish the section with a discussion and a few remarks. Our model extends previous work from random diffusion to nonlocal diffusion, which describes the movement of species from adjacent spatial diffusion to long-distance dispersal. Also, a transient and time-periodic pulse is introduced in the adult to make such structured model describing continuous-discrete process, which increases the reality and universality of application for this model, and the difficulty of mathematical analysis is accordingly raised. Specifically, a juvenile-adult model with nonlocal diffusion and harvesting pulse in moving and heterogeneous environment is proposed to research the dynamics of species in this paper. The global existence and uniqueness of solution to nonlocal problem \eqref{a01} with harvesting pulse in moving boundaries is firstly studied in Theorem \ref{th2.1}. Since the difficulties caused by harvesting pulse and nonlocal operator, we introduce the generalized principal eigenvalue of corresponding eigenvalue problem \eqref{c02} in a fixed domain. The nonincreasing property of the generalized principle eigenvalue related to the length of region is investigated in Lemma \ref{le3.2}(ii). Furthermore, some special conditions for the existence of the principle eigenvalue with positive eigenvalue functions are discovered in Fig. 1 and Remark \ref{rm3.1}. And the strictly decreasing property of the principle eigenvalue related to the length of region and harvesting function is investigated in Lemma \ref{le3.2}(iii). What of particular importance is that the longtime behavior of the solution to problem \eqref{a01} for spreading or vanishing is quite rich, which is caused by dual effects of nonlocal diffusion and harvesting pulse in moving and heterogeneous environment. That is, the sufficient conditions for species to spread or vanish are depicted by $\overline{\lambda}(\infty,H'(0))$ and $\underline{\lambda}(\infty,H'(0))$ in Theorem \ref{th4.9}, and Theorem \ref{th4.10} with some special conditions gives a more concrete dichotomy about spreading-vanishing by $\lambda^*(\infty,H'(0))$. Some sufficient conditions for spreading and vanishing are also given (e.g. $\underline{\lambda}((-h_0,h_0),H'(0))>0$ and small initial value for vanishing in Lemma \ref{le4.7}, $\lambda^*((-h_0,h_0),H'(0))\leq0$ for spreading in Theorem \ref{th4.11}). Spreading-vanishing criteria for expanding capacities are finally
shown in Theorem \ref{th4.11}.

Compared with some other works of mutualistic models, we here consider the longtime behaviors of global solution caused by moving boundaries, harvesting pulse and nonlocal dispersal, which all increase the difficulty of theoretical analysis and diversity of complicated outcomes. In fact, the harvesting rate ($1-H'(0)$) can alter the steady state of solution. For instance, if the principle eigenvalue exists, then for some $H'(0)>0$, we have $\lambda^*(\infty,H'(0))<0$ and species spread (Lemma \ref{le4.8}), while the decreasing property of $H'(0)$ (the increasing property of harvesting rate $1-H'(0)$) may leads to $\lambda^*\geq0$ (Lemma \ref{le3.2}(iii)), which turns species from spread to vanish (Corollary \ref{co4.1}). It indicates that the larger the intensity rate is, individuals goes extinct more likely. Some difficulties caused by nonlocal diffusion and pulse in moving and heterogeneous environment have been overcomed, however, the introduction of the generalized principal eigenvalue leads to a gap in classification of expansion and extinction, which is not conducive to the integrity of the conclusion. Also, whether or not can the extra restrictive conditions about the eigenvalue in several results be removed (see Lemma \ref{le4.4})? And whether does the impact of harvesting timing $\tau$ can be analyzed in this nonlocal and moving boundary problem? All of these are worth further discussion.

\end{document}